\documentclass[11pt]{article}
\usepackage[margin=1.03in]{geometry}
\usepackage{amsmath,amssymb,amsthm,mathtools}
\usepackage{booktabs,array,enumitem,microtype}
\usepackage{longtable}
\usepackage{xcolor}
\usepackage[colorlinks=true,linkcolor=blue!55!black,citecolor=blue!55!black,urlcolor=blue!55!black]{hyperref}
\hypersetup{
  pdftitle={Proof-Carrying Analytic Approximation: Local-to-Global Evidence Transport at Encoding Cost},
  pdfauthor={Andreu Ballus Santacana},
  pdfsubject={Proof-carrying analytic approximation, evidence transport, certified gluing, refinement, computable Banach spaces}
}

\newtheorem{theorem}{Theorem}[section]
\newtheorem{proposition}[theorem]{Proposition}
\newtheorem{lemma}[theorem]{Lemma}
\newtheorem{corollary}[theorem]{Corollary}
\theoremstyle{definition}
\newtheorem{definition}[theorem]{Definition}
\newtheorem{example}[theorem]{Example}
\newtheorem{remark}[theorem]{Remark}

\newcommand{\norm}[1]{\left\lVert #1 \right\rVert}
\newcommand{\abs}[1]{\left\lvert #1 \right\rvert}
\newcommand{\R}{\mathbb{R}}
\newcommand{\N}{\mathbb{N}}
\newcommand{\Q}{\mathbb{Q}}
\newcommand{\Code}{\mathsf{Code}}

\newcommand{\Cert}{\mathbf{Cert}}
\newcommand{\Ev}{\mathsf{Ev}}
\newcommand{\AppCheck}{\mathsf{AppCheck}}
\newcommand{\DistCheck}{\mathsf{DistCheck}}
\newcommand{\PCert}{\mathsf{PCert}}

\newcommand{\aconv}{\operatorname{aconv}}
\newcommand{\range}{\operatorname{ran}}
\newcommand{\id}{\mathrm{id}}

\title{\textbf{Proof-Carrying Analytic Approximation:}\\\textbf{Local-to-Global Evidence Transport at Encoding Cost}}
\author{Andreu Ball\'us Santacana\\
\small Departament de Filosofia, Universitat Aut\`onoma de Barcelona\\
\small \texttt{andreu.ballus@uab.cat}}
\date{}

\begin{document}
\maketitle

\begin{abstract}
Under quasi-uniform refinement, bounded local-encoding hypotheses, and local
$W^{r,2}$ approximation of order $r\ge2$ in a rational piecewise-polynomial presentation
of $W^{1,2}(0,1)$, carrying the complete proof genealogy up to the level required by an
accuracy $\varepsilon$ costs the same asymptotic bit order as the finest-level conventional
coefficient encoding.  If $B_n=\Theta(M_n\beta_n)$ denotes that level-$n$ encoding size,
our compiler transports supplied local approximation and overlap witnesses through exact
partition-of-unity synthesis and geometric refinement to a represented limit with total
certificate size $O(B_{m(\varepsilon)})$, where
$m(\varepsilon)=O(\log(1/\varepsilon)/(r-1))$.  The construction makes no oracle query to
an independently supplied semantic target name ($Q_{\rm target}=0$).  When
$\beta_n=O(n+1)$, this becomes
$O(\varepsilon^{-1/(r-1)}(1+\log(1/\varepsilon)))$.

The surrounding framework is intentionally separated from this resource theorem.  Every
real computable Banach presentation admits a uniformly computable linear isometric
embedding into standard computable $C([0,1])$, with computable inverse on its represented
range.  Complete metric evidence with rational strict slack collapses extensionally to the represented
analytic metric once effective names are available, while chosen evidence transformations
retain construction history and resource information.  For the Lipschitz grammar used
here, qualitative evidence-local lifting is canonical; the nontrivial question is therefore
which evidence is retained and at what cost.
\end{abstract}

\section{Introduction}\label{sec:intro}

Proof-carrying analytic approximation separates two questions that ordinary approximation
usually identifies.  The first is extensional: how an analytic object is represented and
approximated by finite data.  The second is genealogical: which finite derivation produced
an approximant, how local evidence was transported into global evidence, and what it costs
to retain that history through refinement.  Standard represented analysis answers the
first question.  This paper isolates a concrete mathematical form of the second.  The
central problem is whether standard analytic inequalities can be compiled into finite proof
objects without discarding incoming evidence and recertifying the target from scratch.

The distinction is necessary because generic certificate existence is cheap once a full
Type-2 name is already available.  If a true metric inequality has strict rational slack,
a checker allowed to read both names can search to a finite stage at which the Cauchy tails
and a rational enclosure certify it.  Under that completeness condition the separated
evidence metric therefore collapses to the represented analytic metric.  What can remain
proof-relevant is not a new extensional point, but construction history, selected
derivations, proof size, verification work, source use, and the choice of evidence
transformer.  This is close in spirit to assembly semantics: ordinarily one retains the
underlying function for which some tracker exists; here a selected evidence transformer is
kept as additional arrow data.  Throughout, the proof language and its binary encoding are
part of the presentation, and no proof-size comparison is treated as representation
invariant without an explicit translation theorem.

The analytic representation layer supplies a uniform ambient language.  Theorem~\ref{thm:U}
and Corollary~\ref{cor:finite-linear}, together called UELAT below, give a computable linear
isometric embedding of every real computable Banach presentation into the standard
computable $C([0,1])$, with computable inverse on the represented range, and hence a common
finite rational polygonal language for approximating embedded points.  This theorem is
standalone: the quantitative Sobolev construction later works in its native rational
piecewise-polynomial presentation.  The role of UELAT is instead to show that the finite-code
requirement of the certificate framework is non-vacuous uniformly across real computable
Banach presentations.

The main local-to-global result is more concrete.  In $W^{1,2}(0,1)$, rational
piecewise-polynomial codes make finite $H^1$ arithmetic exact.  Supplied local approximation
and overlap witnesses are compiled through rational partition-of-unity synthesis and a
geometric refinement sequence into a shared proof DAG.  The output at precision
$\varepsilon$ uses no independently supplied target name: $Q_{\rm target}=0$ means zero
oracle calls to such a semantic target name, not zero use of supplied source evidence.
If $B_n=\Theta(M_n\beta_n)$ is the declared conventional coefficient-wise encoding size at
level $n$, the complete genealogy through level $m(\varepsilon)$ has size
$O(B_{m(\varepsilon)})$.  The mechanism is the geometric collapse
\[
 \sum_{j\le n}M_j\beta_j=O(M_n\beta_n):
\]
coarser refinement levels sum to a constant multiple of the finest one.  The
$h^{r-1}$ scale law determines the required level, bounded same-level incidence controls
local dependencies, and the encoding hypothesis prevents proof payload from outgrowing the
ordinary approximation data.  This is an upper-bound theorem relative to a fixed encoding;
no information-theoretic lower bound asserting that target queries or this certificate size
are necessary is claimed.

Exact finite arithmetic and strict categorical bookkeeping play different roles.  Exact
rational synthesis makes code-level reconstruction and compatibility equations literal.
Strict identity and associativity of evidence composition are supplied separately by a
normalized endpoint-indexed proof-spine convention: reflexivity is the empty spine and
triangle composition is flattened concatenation.  Chosen transport acts on these normalized
spines, so the later $\mathbf{Cat}$-valued organization is strict rather than merely strict
up to syntactic rebracketing.

Contextual Choice supplies an operational admissibility discipline around these
constructions.  A declared certified development starts from primitive certified objects
and specified evidence-transforming constructors; local selections are promoted globally
only through proved reconstruction, gluing, or refinement interfaces with their required
moduli.  The least generated closure then supports invariant-preservation and internal
exclusion arguments, with Corollary~\ref{cor:BT} as a concrete finite-measure example.
Nothing in this discipline forbids an external classical choice outside the declared
signature.  Exact compatible sections, after extensional collapse, retain their ordinary
sheaf semantics.

\paragraph{Contributions.}
The paper makes four claims, separated by mathematical layer.
\begin{enumerate}[label=(\roman*),leftmargin=2.2em,itemsep=2pt,topsep=3pt]
\item A uniform effective linear Banach--Mazur representation theorem into standard
computable $C([0,1])$, with computable inverse on the represented range, and the resulting
common finite rational code language for embedded points.
\item A qualitative boundary theorem for finite evidence: generic complete strict-slack
evidence collapses extensionally, while chosen lifts retain proof-relevant data and the
finite-code Lipschitz grammar admits canonical evidence-local transport.
\item A quantitative rational $W^{1,2}(0,1)$ compiler that carries supplied local evidence
through exact synthesis and refinement to a represented limit with $Q_{\rm target}=0$ and
complete provenance size $O(B_{m(\varepsilon)})$.  Under the standard linear bit schedule
$\beta_n=O(n+1)$ this yields
$O(\varepsilon^{-1/(r-1)}(1+\log(1/\varepsilon)))$.
\item An operational Contextual Choice discipline of generated genealogical closure under
declared constructors, together with invariant preservation and the internal boundary
example of Corollary~\ref{cor:BT}; exact sheaf gluing remains extensional.
\end{enumerate}

\paragraph{Organization.}
Section~\ref{sec:evidence} fixes represented objects, finite evidence, the access-sensitive
collapse result, and the normalized proof-spine convention.  Section~\ref{sec:universal}
proves the universal effective linear representation theorem.  Section~\ref{sec:lifts}
develops chosen proof-relevant transport and its qualitative saturation boundary.
Sections~\ref{sec:reconstruction}--\ref{sec:descent} give the rational $W^{1,2}$
reconstruction and refinement compilers and prove the same-order resource theorem.
Section~\ref{sec:ccp} gives the generated-universe semantics; the following boundary and
related-work sections place the result against standard computation, certification, and
sheaf theory.  The public-version correspondence is consolidated in
Section~\ref{sec:versions}; no theorem depends on that history.

\section{Represented analytic objects and finite evidence}\label{sec:evidence}

\subsection{The represented-space baseline}

\begin{definition}[Real computable Banach presentation]\label{def:cban}
A \emph{real computable Banach presentation} (real meaning scalar field $\R$) is a Banach
space $X$ together with a dense sequence $(e_i)_{i\in\N}$ such that rational linear
combinations of the $e_i$ form an effectively enumerable rational vector space, addition
and rational scalar multiplication on that core are computable, and the norm of a
rational-core vector is a computable real, uniformly in its code.  A standard name of
$x\in X$ is a fast Cauchy sequence $(p_n)$ of rational-core codes satisfying
\[
 \norm{x-\rho(p_n)}\le 2^{-n}.
\]
Computability of points and maps is understood through this Cauchy representation in the
usual Type-2 sense.  When a theorem below is uniform in the presentation itself, its input
is an index for the enumerated core operations and norm procedure just specified; correctness
is promised on indices that actually present Banach spaces in this sense.
\end{definition}

This is the standard computable-metric/computable-Banach starting point; nothing in the
certificate layer changes the represented identity of a point~\cite{PourElRichards,Weihrauch}.

\begin{definition}[Certificate enrichment]\label{def:enrichment}
Let $(X,\delta_X)$ be a represented normed space and let $D_X\subseteq X$ be the
represented subdomain under discussion.  A \emph{certificate enrichment}
$\mathfrak X$ consists of:
\begin{enumerate}[label=(\roman*)]
\item an effectively enumerable finite-code set $\Code_X$ with decoder
$\rho_X:\Code_X\to D_X$ whose image is dense in $D_X$;
\item a terminating approximation checker $\AppCheck_X$ whose acceptance
\[
 \AppCheck_X(\nu,p,q,V)=1
\]
implies $\norm{\delta_X(\nu)-\rho_X(p)}\le q$;
\item a terminating distance checker $\DistCheck_X$ whose acceptance
\[
 \DistCheck_X(\nu,\mu,q,W)=1
\]
implies $\norm{\delta_X(\nu)-\delta_X(\mu)}\le q$;
\item a declared typed finite-evidence signature $\Sigma_X$ containing the structural
constructors used in the development (identity, symmetry, weakening, triangle, and any
declared bounded-operator, finite-sum, overlap, or exact-arithmetic rules).
\end{enumerate}
Clauses (ii)--(iii) are the \emph{soundness} requirement.  Structural distance witnesses
are stored in a fixed normalized endpoint-indexed spine form: reflexivity is the empty
spine and triangle composition is flattened concatenation, with announced rational bounds
added.  Declared transport constructors preserve the empty spine and concatenation.  Thus
the witness composition used below is unital and associative as literal morphism equality,
not merely up to syntactic rebracketing.  Primitive-step validity remains checker-mediated.
The proof language, normalization convention, and binary encoding are all part of the
presentation.
\end{definition}

\begin{definition}[Certificate and certificate system]\label{def:certificate}
For a name $\nu$ and rational $\varepsilon>0$, a certificate is a finite record
\[
  (p,\bar\varepsilon,V),\qquad
  p\in\Code_X,\quad 0\le\bar\varepsilon<\varepsilon,
\]
with $\AppCheck_X(\nu,p,\bar\varepsilon,V)=1$.  Under a fixed binary encoding its length is
\[
 |(p,\bar\varepsilon,V)|
 :=|p|+\ell_{\Q}(\bar\varepsilon)+|V|.
\]
Here $\ell_{\Q}$ is the bit length of the chosen binary rational encoding.
A certificate system over $\nu$ is a uniform terminating procedure which returns such a
record for every rational $\varepsilon>0$.  Unless a later access discipline says otherwise,
the system may read its supplied source name $\nu$; target-name use is accounted for
separately in Definitions~\ref{def:locality}--\ref{def:resources}.
\end{definition}

The finite record is intentionally stronger than an approximant alone: it includes the
finite reason the declared checker accepts the error claim.  It is nevertheless weaker
than a new representation of the analytic point, because the represented name remains the
extensional carrier.

\subsection{Canonical strict-slack evidence and collapse}

\begin{definition}[Strict-slack completeness]\label{def:slackcomplete}
A certificate enrichment is \emph{strict-slack complete} if, uniformly in names $\nu,\mu$
and rational $q>\norm{\delta(\nu)-\delta(\mu)}$, one can produce finite accepted distance
evidence at bound $q$, and similarly, for every code $p$ and rational
$q>\norm{\delta(\nu)-\rho(p)}$, one can produce accepted approximation evidence at bound
$q$.
\end{definition}

\begin{proposition}[Canonical slack certification]\label{prop:slack}
For a computable metric presentation (in particular, a computable Banach presentation) whose checker is allowed to read the supplied
Cauchy names, strict-slack completeness can be realized by finite precision.  In
particular, if $d(\delta\nu,\delta\mu)<q$, then there is a finite stage $n$ at which a
computable rational upper enclosure of the distance between the $n$th approximants,
together with the Cauchy tails, lies strictly below $q$; storing $n$ and the corresponding
finite arithmetic data is a terminating witness.
\end{proposition}

\begin{proof}
Let $d=d(\delta\nu,\delta\mu)<q$ and put $\gamma=q-d>0$.  Fast Cauchy names give
approximants $x_n,y_n$ with $d(x_n,\delta\nu),d(y_n,\delta\mu)\le2^{-n}$.  By two triangle inequalities, for every $n$,
\[
 d(\delta\nu,\delta\mu)\le d(x_n,y_n)+2^{1-n},
\]
and the right-hand side converges to $d(\delta\nu,\delta\mu)$.
The finite-stage distance is a computable real, so we can compute rational upper
enclosures converging to it.  For sufficiently large $n$ one of these upper enclosures
plus $2^{1-n}$ is $<q$.  Search over $n$ and over enclosure precision until this strict
rational inequality is observed.  The search terminates because $\gamma>0$.  The
approximation case is identical with one Cauchy tail.  This is the standard finite-slack
mechanism underlying complete metric checkers.
\end{proof}

The proposition is deliberately conditional on effective access to the names.  A verifier
restricted to a coarser information interface need not be complete.  The corresponding
lower-bound question -- whether target queries are necessary when no construction genealogy
is supplied -- is not addressed here; finite-query indistinguishability is studied in
information-based complexity and the SCI literature~\cite{TWW,ColbrookHansen,Sorg2026}.

\begin{definition}[Proof-relevant evidence category]\label{def:evcat}
For a certificate enrichment $\mathfrak X$, let $\Cert_{\rm ev}(\mathfrak X)$ have
objects $C=(\nu,c)$ consisting of a represented name and a certificate system over it.
A morphism
\[
 C=(\nu,c)\longrightarrow D=(\mu,d)
\]
is a pair $(q,W)$ with $q\in\Q_{\ge0}$ and an accepted normalized endpoint-indexed
witness spine $W$ from $\nu$ to $\mu$ whose certified intrinsic bound is at most $q$.
The identity is the empty spine at announced bound $0$; composition concatenates spines
and adds announced bounds.  By the normalization convention of
Definition~\ref{def:enrichment}, these operations satisfy the category laws strictly.
The certificate systems $c,d$ are intentionally inert in distance-morphism acceptance;
they remain object data and are acted on by the chosen evidence lifts of Section~\ref{sec:lifts}.
\end{definition}

\begin{definition}[Evidence pseudometric]\label{def:evmetric}
For evidence objects $C,D$, define
\[
 d_{\Cert}(C,D)
 :=\inf\{q\in\Q_{\ge0}:\exists W\;(q,W):C\to D\},
\]
with $\inf\varnothing=\infty$.  The forgetful map $U_X$ sends $(\nu,c)$ to
$\delta_X(\nu)$.  The symmetry constructor makes this particular Lawvere-style
pseudometric symmetric.  The infimum-over-weighted-arrows construction belongs to the
Lawvere generalized-metric tradition~\cite{Lawvere73}; no new enriched universal
property is claimed here.
\end{definition}

\begin{theorem}[Conditional extensional collapse]\label{thm:collapse}
If $\mathfrak X$ is sound and strict-slack complete for distance claims (as in Proposition~\ref{prop:slack} when full name access is available), then
\[
 d_{\Cert}(C,D)=\norm{U_X(C)-U_X(D)}
\]
for all evidence objects $C,D$.  Consequently the zero-distance quotient of the evidence
pseudometric is canonically isometric to the represented points that admit certificate
systems, equipped with the ambient metric.
\end{theorem}

\begin{proof}
Soundness gives $\norm{U_X(C)-U_X(D)}\le d_{\Cert}(C,D)$.  If
$r=\norm{U_X(C)-U_X(D)}$ and $q>r$ is rational, strict-slack completeness provides an
accepted morphism at bound $q$, so $d_{\Cert}(C,D)\le q$.  Taking the infimum over
rational $q>r$ gives equality.  Therefore zero evidence distance is exactly equality of
the underlying represented points, and the quotient statement follows.
\end{proof}

\begin{remark}[What collapses]
The theorem does not make finite evidence meaningless.  It says that once effective
metric access is already supplied, unrestricted qualitative strict-slack certificates do
not change the \emph{extensional carrier}.  Construction history, selected derivations,
proof size, verification work, source use, and the choice of evidence transformer can
still differ above the same analytic point.
\end{remark}

\section{Universal effective linear representation}\label{sec:universal}

The ambient universal representation component is developed in this section.  We reserve
the acronym UELAT -- Universal Embedding and Linear Approximation -- for
Theorem~\ref{thm:U} together with its finite-code corollary,
Corollary~\ref{cor:finite-linear}.  Here ``linear approximation'' means finite linear
representation of points in the ambient space, not finite-rank approximation of the
identity; see Remark~3.5.  The theorem is a uniform computable linear isometric embedding
into one standard ambient Banach space; the code language follows from that embedding.
This definition is intrinsic to the present paper.

The classical Banach--Mazur theorem states that every separable Banach space embeds linearly
and isometrically into $C([0,1])$; see, for example, Albiac and Kalton~\cite{AlbiacKalton}.
Metric computable universality of the standard computable $C([0,1])$ is also
known~\cite{BagavievEtAl,DowneyMelnikov}.  The issue for a linear computable embedding is
that the usual norming-functional proof invokes exact Hahn--Banach selection, which is not
computable in general~\cite{Brattka2008,BrattkaSorg2026}.  We avoid such a selector.

\begin{lemma}[Effective approximate Hahn--Banach]\label{lem:epsHB}
Let $X$ be a real computable Banach presentation, let $v$ be a nonzero computable vector,
and let $\eta>0$ be rational.  Uniformly on this promised domain one can compute a Type-2
name of a bounded linear functional $g:X\to\R$ such that
\[
 g(v)=\norm v,
 \qquad \norm g\le 1+\eta.
\]
\end{lemma}

\begin{proof}
Because the norm is computable and $v\ne0$, search for $s$ until $2^{-s}<\norm v$ is
certified.  This search is uniform and terminates on the stated domain.  Let $Y=\R v$.
The map $a\mapsto av$ gives $Y$ a computable one-dimensional presentation; the certified
lower bound also makes the inverse coordinate map effective, since
\[
 |a-b|\le 2^s\norm{av-bv}.
\]
Hence
\[
 f_0:Y\to\R,\qquad f_0(av)=a\norm v,
\]
is a computable bounded functional with $\norm{f_0}=1$ and located kernel $\{0\}$, the locatedness hypothesis required by the constructive extension theorem used next.

The constructive $\varepsilon$-Hahn--Banach theorem extends such a functional with an
arbitrarily prescribed norm loss; see~\cite[Chap.~7, Thm.~4.6, p.~342]{BishopBridges}.
Brattka explicitly notes that Bishop's fully constructive $\varepsilon$-version, in which
norm preservation is relaxed to $\norm g\le\norm{f_0}+\eta$, transfers to the represented
computable setting~\cite{Brattka2005}.  Applying that effective construction to the
computable presentation of $X$, the computable inclusion $Y\hookrightarrow X$, $f_0$, and
the rational parameter $\eta$ produces a Type-2 name of an extension $g$.  Extension gives
$g(v)=f_0(v)=\norm v$, and the norm estimate gives $\norm g\le1+\eta$.

All operations just described -- the positive-norm search, the one-dimensional presentation,
the functional $f_0$, and the effective approximate extension -- are uniform in the
presentation index of Definition~\ref{def:cban}, a name of $v$, and $\eta$.  No exact norm-preserving Hahn--Banach
selector is computed.
\end{proof}

\begin{theorem}[Effective linear Banach--Mazur universality]\label{thm:U}
There is a single effective procedure which, from an index promised to be a valid real
computable Banach presentation $X$ in the sense of Definition~\ref{def:cban}, produces a Type-2 computable linear isometric embedding
\[
 J_X:X\longrightarrow C([0,1])
\]
into the standard computable presentation of $C([0,1])$.  The inverse
\[
 J_X^{-1}:J_X[X]\longrightarrow X
\]
is computable on the represented range.
\end{theorem}

\begin{proof}
We separate the construction into seven uniform steps.

\smallskip
\noindent\emph{Step 1: a computable $1$-norming family.}
Enumerate the rational linear span of the fundamental sequence as $(v_i)_{i\in\N}$.
For pairs $(i,s)$ semidecide $\norm{v_i}>2^{-s}$; every nonzero rational-core vector is
eventually listed with a certified positive lower bound.  For such $v=v_i$ and
$\eta_k=2^{-k}$, Lemma~\ref{lem:epsHB} gives a computable extension $g_{v,k}$ satisfying
\[
 g_{v,k}(v)=\norm v,
 \qquad \norm{g_{v,k}}\le1+\eta_k.
\]
Put
\[
 u_{v,k}:=\frac{g_{v,k}}{1+\eta_k}.
\]
Then $\norm{u_{v,k}}\le1$ and
\[
 u_{v,k}(v)=\frac{\norm v}{1+\eta_k}.
\]
Flatten these functionals into a uniformly computable sequence $(u_j)$ in the dual unit
ball.  It is $1$-norming:
\[
 \norm x=\sup_j\abs{u_j(x)}.
\]
Indeed the upper bound follows from $\norm{u_j}\le1$.  For $x\ne0$ and $\varepsilon>0$,
choose a rational-core $v$ sufficiently close to $x$ and then $k$ sufficiently large.
Since $u_{v,k}$ has norm at most one,
\[
 u_{v,k}(x)
 \ge \frac{\norm v}{1+\eta_k}-\norm{x-v},
\]
and the right-hand side approaches $\norm x$ as $v\to x$ and $k\to\infty$.

\smallskip
\noindent\emph{Step 2: a coordinate compactum for the dual ball.}
Choose positive rational numbers $c_i>\max(1,\norm{v_i})$ effectively and set
$z_i=v_i/c_i$, so $\norm{z_i}<1$ and the rational span of $(z_i)$ is dense in $X$.
Inside the computable Hilbert cube $H=[-1,1]^{\N}$ define
\[
 K_X:=\left\{a=(a_i):
 \abs{\sum_i q_i a_i}\le
 \norm{\sum_i q_i z_i}
 \text{ for every finitely supported }(q_i)\in\Q^{(\N)}
 \right\}.
\]
Each $a\in K_X$ defines a rational-linear functional on
$\operatorname{span}_{\Q}\{z_i\}$ by $\sum q_i z_i\mapsto\sum q_i a_i$.  This is
well defined even though the $z_i$ need not be linearly independent: if two rational
combinations represent the same vector, applying the defining inequality to their
difference forces the corresponding scalar difference to be $0$.  The displayed inequality
then gives norm at most one, so the functional extends uniquely by continuity to a member
of $B_{X^*}$.  Conversely every $\phi\in B_{X^*}$ gives coordinates
$a_i=\phi(z_i)$ and belongs to $K_X$.  Thus $K_X$ is precisely the weak-* dual unit ball
in these coordinates.  Failure of one defining inequality is effectively open; hence
$K_X$ is uniformly co-c.e. closed in $H$.  It is nonempty because it contains the zero
sequence.  This supplies the negative compact-set information used below; the positive
dense information needed for a computable metric presentation is supplied in Step~3.

\smallskip
\noindent\emph{Step 3: effective dense points in $K_X$.}
Enumerate all finite rational absolutely convex combinations
\[
 w=\sum_{j\le N} r_j u_j,
 \qquad r_j\in\Q,\quad \sum_j\abs{r_j}\le1.
\]
Their coordinate sequences $(w(z_i))_i$ are uniformly computable points of $K_X$.
Because $(u_j)$ is $1$-norming, its polar in $X$ is the closed unit ball.  The bipolar
theorem therefore gives
\[
 \overline{\aconv\{u_j:j\in\N\}}^{\;w^*}=B_{X^*}.
\]
Hence the enumerated computable points are weak-* dense in $K_X$.  On the bounded dual
ball, the inherited Hilbert-cube coordinate metric metrizes the weak-* topology because the
rational span of the $z_i$ is dense in $X$; therefore these points are also dense for the
metric used in Step~4.

\smallskip
\noindent\emph{Step 4: computable compactness.}
The positive information from Step~3 is essential here: co-c.e.\ closedness by itself does
not supply computable points, and a nonempty co-c.e.\ compact set can fail to contain any
computable point.  Equip $K_X$ with the inherited computable Hilbert-cube metric and the
dense sequence from Step~3.  To find a finite $2^{-n}$-net, enumerate finite families of dense points until
their $2^{-n}$-balls cover $K_X$.  Coverage is semidecidable: if $U$ is the proposed
finite union, then
\[
 K_X\subseteq U
 \quad\Longleftrightarrow\quad
 H\subseteq (H\setminus K_X)\cup U,
\]
and the right-hand side is a c.e.-open cover of the computably compact Hilbert cube.
Compactness and density guarantee termination.  Thus $K_X$ is a nonempty computably
compact computable metric space.

\smallskip
\noindent\emph{Step 5: evaluate over the whole dual compactum.}
By the effective Hausdorff--Alexandroff theorem for computably compact metric spaces,
there is a computable continuous surjection
\[
 q_X:2^{\N}\twoheadrightarrow K_X
\]
uniformly from the computably compact presentation~\cite[Prop.~4.1]{BrattkaDeBrechtPauly};
see also~\cite[Thm.~3.40]{DowneyMelnikov}.
For $x\in X$ define
\[
 (T_Xx)(p):=\phi_{q_X(p)}(x),
 \qquad p\in2^{\N}.
\]
Evaluation is uniformly computable.  Given $x$, approximate it by a finite rational
combination $z=\sum_i q_i z_i$; the coordinates of $q_X(p)$ compute
$\phi_{q_X(p)}(z)$, while the unit-ball bound controls the residual by
$\norm{x-z}$.  Currying gives a computable map
\[
 T_X:X\to C(2^{\N}).
\]
It is linear.  Since every value of $q_X$ lies in the dual unit ball,
$\norm{T_Xx}_{\infty}\le\norm x$.  Conversely each $u_j$ from Step~1 belongs to
$B_{X^*}=\range(q_X)$, so
\[
 \norm{T_Xx}_{\infty}
 \ge \sup_j\abs{u_j(x)}
 =\norm x.
\]
Thus $T_X$ is isometric without selecting, or invoking in this verification, an exact
Hahn--Banach norming functional.

\smallskip
\noindent\emph{Step 6: from Cantor space to $[0,1]$.}
Identify $2^{\N}$ computably with the ternary Cantor set $C\subset[0,1]$.  This gives a
computable linear isometry $C(2^{\N})\cong C(C)$.  For $f\in C(C)$ define $Ef\in C([0,1])$
by $Ef=f$ on $C$ and by affine interpolation on each complementary interval of $C$.
The operator $E$ is linear and
\[
 \norm{Ef}_{\infty}=\norm f_{\infty},
\]
because every gap value is a convex combination of endpoint values.  It is computable in
the standard representations: an effective modulus of uniform continuity controls the
small gaps, while finitely many larger gaps are handled by rational polygonal
approximation.  Set
\[
 J_X:=E\circ T_X.
\]
Then $J_X$ is computable, linear, and isometric.

\smallskip
\noindent\emph{Step 7: computable inverse on the range.}
The images $J_X(v_i)$ of the rational core form a computable dense sequence in $J_X[X]$.
Given a name of $y$ promised to lie in the represented range, say $y=J_Xx$, and $n$, search for $i$ such that
\[
 \norm{J_X(v_i)-y}_{\infty}<2^{-n}.
\]
Strict inequality of the computable $C([0,1])$ distance is semidecidable, and density
ensures termination.  By isometry,
\[
 \norm{v_i-x}<2^{-n}.
\]
The selected rational-core vectors form a computable Cauchy name of $x$, proving that
$J_X^{-1}$ is computable on its represented range.
\end{proof}

\begin{remark}[Exact Hahn--Banach is not being computed]
The construction is consistent with the noncomputability of exact norm-preserving
Hahn--Banach extension.  Approximate extensions are used only to produce a countable
$1$-norming family.  That family gives exact norm recovery; the effectively presented dual
unit ball supplies the compact domain needed to land in $C(2^{\N})$ and then $C([0,1])$.  Exact Hahn--Banach retains nontrivial Weihrauch complexity even
on fixed classical Banach spaces~\cite{BrattkaSorg2026}.
\end{remark}

\begin{corollary}[Universal finite linear representation]\label{cor:finite-linear}
Every computable point $x$ of every real computable Banach presentation can be transported
computably to $J_Xx\in C([0,1])$ and approximated to arbitrary precision by a finite
rational polygonal (equivalently, rational hat-function) code.  The individual polygonal approximants need not lie in $J_X[X]$; their represented
limit does, and pulling that limit back through $J_X^{-1}$ recovers $x$.
\end{corollary}

\begin{remark}[Linear approximation means point representation]
The corollary concerns approxi\-mation of individual embedded points by finite linear
combinations from a dense universal code family.  It does \emph{not} assert that every
Banach space has the approximation property, i.e.\ that its identity is uniformly
approximable on compact subsets by finite-rank operators, still less the bounded
approximation property.  Enflo's counterexample shows that the approximation property
fails for some Banach spaces~\cite{Enflo}.
\end{remark}

\begin{remark}[Role of the universal representation theorem]\label{rem:role-universal}
Definition~\ref{def:enrichment}(i) posits an effectively enumerable dense finite-code
language.  Theorem~\ref{thm:U} and Corollary~\ref{cor:finite-linear} construct such a
language uniformly for every real computable Banach presentation, inside one standard
ambient space; in that sense the theorem is a non-vacuity result for the general certificate
framework.  The quantitative Sobolev compiler in
Sections~\ref{sec:reconstruction}--\ref{sec:descent} does not factor its PUFEM estimates
through $J_X$ and instead uses its native exact rational presentation.  Transporting an
arbitrary enrichment along $J_X$ would additionally require finite-code compilers between
the chosen presentations; no general quantitative transport theorem of that kind is asserted here.
\end{remark}

\section{Proof-relevant evidence transport}\label{sec:lifts}

\subsection{Finite-code realizable Lipschitz maps}

\begin{definition}[Finite-code realizable Lipschitz map]\label{def:code-realizable}
Let $\mathfrak X,\mathfrak Y$ enrich represented normed spaces $X,Y$.  A map on the
underlying spaces, $T:X\to Y$, is \emph{finite-code realizable Lipschitz} when the
enrichments supply:
\begin{enumerate}[label=(\roman*)]
\item a rational $\Lambda_T\ge0$ and a stored derivation of
$\norm{Tx-Ty}\le\Lambda_T\norm{x-y}$;
\item a Type-2 realizer $T^\sharp$ on represented names;
\item a uniform finite-code compiler which, given $p\in\Code_X$ and rational $\eta>0$,
returns a pair $(q,E_T(p,\eta))$ with $q=\tau_T(p,\eta)\in\Code_Y$ and finite evidence
for
\[
 \norm{T\rho_X(p)-\rho_Y(q)}\le\eta.
\]
\end{enumerate}
The name realizer and finite-code compiler are coherent through the same underlying map
$T$; no stronger syntactic relation is required.  The compiler is \emph{exact} when
$\rho_Y(\tau_T(p,\eta))=T\rho_X(p)$ for all $p$, in which case $\eta$ is inert and
the defect evidence may be trivial.  When quantitative use is intended, a size bound for
$(q,E_T)$ is part of the declared compiler data.  Sections~\ref{sec:reconstruction}--\ref{sec:descent}
derive such bounds directly for their concrete compilers rather than assuming a general
resource-composition theorem at this level.
\end{definition}

\begin{theorem}[Qualitative local-transport saturation]\label{thm:saturation}
Assume the evidence grammar admits tagged source subproofs, the stored Lipschitz rule,
triangle, and weakening.  Every finite-code realizable Lipschitz map has a canonical
evidence-local lift
\[
 T_*:\Cert_{\rm ev}(\mathfrak X)\longrightarrow\Cert_{\rm ev}(\mathfrak Y).
\]
At target tolerance $\varepsilon$, one source approximation certificate at
\[
 \alpha_T(\varepsilon)=\frac{\varepsilon}{3\max(1,\Lambda_T)}
\]
and one finite-code compiler call at defect $\varepsilon/3$ suffice.  The lift makes no semantic query to the target point beyond the output of the declared
code compiler; this is the local form of the $Q_{\rm target}=0$ audit in
Theorem~\ref{thm:orderneutral}.
\end{theorem}

\begin{proof}
Suppose source evidence $V$ proves
$\norm{x-\rho_X(p)}\le r$.  The stored Lipschitz rule gives
\[
 \norm{Tx-T\rho_X(p)}\le\Lambda_T r.
\]
The finite-code compiler supplies $E_T$ proving
$\norm{T\rho_X(p)-\rho_Y(q)}\le\eta$.  One triangle node produces
\[
 \norm{Tx-\rho_Y(q)}\le\Lambda_T r+\eta.
\]
This defines the approximation transformer
\[
 \Xi_T(V,E_T)=\operatorname{triangle}(\operatorname{Lip}_T(V),E_T).
\]
Likewise source distance evidence $W$ at bound $r$ is sent to
$\Theta_T(W)=\operatorname{Lip}_T(W)$ at bound $\Lambda_T r$.  Taking the stated $\alpha_T(\varepsilon)$ and $\eta=\varepsilon/3$ gives at most
$2\varepsilon/3<\varepsilon$; the three-way split is only a convenient strict-slack
normalization.  On distance morphisms the transport acts on normalized witness spines,
preserving the empty spine and concatenation by Definition~\ref{def:enrichment}; hence
identity and composition hold strictly.
\end{proof}

\begin{remark}[Consequence]
Within this Lipschitz calculus, an additional axiom asserting the mere existence of
approximation- and distance-evidence transformers is redundant once the analytic estimate
and finite-code action are supplied.  Thus evidence locality by itself does not create a
smaller qualitative morphism class here.  A nontrivial distinction must concern which
lift is chosen, which incoming evidence it retains, or the resources used to transport it.
\end{remark}

\subsection{Chosen lifts and non-faithfulness}

\begin{definition}[Analytic maps with chosen lifts]\label{def:chosen}
Let $\mathbf{CAn}$ have certificate enrichments as objects and finite-code realizable
Lipschitz maps as morphisms.  Let $\mathbf{CAn}^{\uparrow}$ have the same objects and
morphisms pairs $(T,\widetilde T)$, where $T$ is such a map and $\widetilde T$ is a
chosen evidence functor satisfying
\[
 U_Y\widetilde T=T U_X.
\]
Composition is componentwise.  The forgetful functor
\[
 P:\mathbf{CAn}^{\uparrow}\to\mathbf{CAn}
\]
is the identity on objects and sends $(T,\widetilde T)$ to $T$.
\end{definition}

\begin{definition}[Evidence regularity]\label{def:regular}
A certificate enrichment is \emph{evidence regular} if each finite code $p$ has an exact
self-certificate over its canonical constant name and accepted approximation evidence can
be promoted to distance evidence between the represented point and that canonical name.
\end{definition}

\begin{proposition}[Non-faithfulness of analytic forgetting]\label{thm:nonfaithful}
Let $\mathfrak X$ be a nonempty evidence-regular enrichment.  Then $P$ is not faithful.
More strongly, there is a monoid of pairwise distinct endofunctor lifts
$R_k:\Cert_{\rm ev}(\mathfrak X)\to\Cert_{\rm ev}(\mathfrak X)$, $k\in\N$, satisfying
\[
 U_XR_k=U_X,\qquad R_kR_\ell=R_{k+\ell},
\]
all lying over the analytic identity.
\end{proposition}

\begin{proof}
For a certificate system $c$, define
\[
 (R_kc)(\varepsilon):=c(2^{-k}\varepsilon).
\]
If $c(2^{-k}\varepsilon)=(p,\bar\varepsilon,V)$, then
$\bar\varepsilon<2^{-k}\varepsilon\le\varepsilon$, so the same finite record is a valid
certificate at tolerance $\varepsilon$.  On distance morphisms let $R_k$ act as the
identity.  This is well defined because morphism acceptance depends on the represented
names and the supplied distance witness, not on the schedule by which the object-level
certificate systems answer tolerances.  The displayed identities are immediate.

To prove distinctness, choose a finite code $p$ and its exact self-witness $V_p^0$.
Weakening yields a certificate system
\[
 c_*(\varepsilon)=\left(p,\frac{\varepsilon}{2},
       \operatorname{weak}(V_p^0,\varepsilon/2)\right).
\]
Then
\[
 (R_kc_*)(\varepsilon)
 =\left(p,\frac{\varepsilon}{2^{k+1}},\ldots\right),
\]
so $R_k\ne R_\ell$ for $k\ne\ell$.  In particular
$(\id_X,R_0)\ne(\id_X,R_1)$ but both map under $P$ to $\id_X$.
\end{proof}

This proposition is structural rather than a headline complexity separation: forgetting data
that were deliberately made part of an arrow is expected to lose intensional information.
Its role is to identify exactly what survives the extensional collapse.  In ordinary
assembly semantics, a morphism is the underlying function for which some tracker exists;
the chosen tracker is not itself arrow identity~\cite{vanOosten}.  The selected
$\widetilde T$ in $\mathbf{CAn}^{\uparrow}$ is therefore additional proof-relevant data.

\begin{proposition}[Grothendieck organization]\label{prop:opfibration}
The strict assignment
\[
 \Ev:\mathbf{CAn}^{\uparrow}\to\mathbf{Cat},\qquad
 \mathfrak X\mapsto\Cert_{\rm ev}(\mathfrak X),\quad
 (T,\widetilde T)\mapsto\widetilde T,
\]
has the usual Grothendieck construction, whose projection to
$\mathbf{CAn}^{\uparrow}$ is a split opfibration.
\end{proposition}

\begin{proof}
Strictness is by construction: composition in $\mathbf{CAn}^{\uparrow}$ is componentwise and each chosen lift is a genuine functor on the normalized evidence categories.  The result is the standard Grothendieck construction for a strict covariant $\mathbf{Cat}$-valued functor.  The chosen opcartesian lift of
$(T,\widetilde T)$ at $c$ is $((T,\widetilde T),\id_{\widetilde T(c)})$.
\end{proof}

\subsection{Evidence locality and resource profiles}

\begin{definition}[Evidence-local transformer]\label{def:locality}
An evidence transformer is \emph{evidence local} when an output certificate at tolerance
$\varepsilon$ is generated from finitely many supplied source certificates, finite codes,
finite operation data, and explicitly declared source-name queries, without invoking an
independent generic certificate generator on the semantic target name.  Evidence locality
is an access restriction; it is not a claim of novelty over tracker semantics or
finite-information computation.
\end{definition}

\begin{definition}[Resource profile]\label{def:resources}
For a fixed presentation, proof language, binary encoding, and oracle-machine model, we
distinguish
\[
 (G,V,T,S,Q;L).
\]
Here $G$ is certificate-generation cost, $V$ verification cost, $T$ evidence-transport
cost, $S$ encoded output proof size, and $Q$ represented-name/oracle use.  A query means
one oracle call requesting a represented name at a specified precision; the zero-query
statements below are therefore independent of finer conventions for charging nonzero
lookahead.  The symbol $L$ after the semicolon is a predicate recording the evidence-local
access discipline, not a numerical resource.  The concrete compiler results below
explicitly bound $S$, $V$, and target-query use; the remaining components provide context
for the declared machine model.  All quantitative claims are representation and encoding
dependent.  Resource-monoid semantics was introduced by Dal Lago and Hofmann~\cite{DalLagoHofmann};
Hofmann and Ledent later use the framework to majorize normal-form size rather than only
execution time~\cite{HofmannLedent}.  The vector above is only bookkeeping for the concrete
analytic compiler, not a proposed replacement for those semantics.
\end{definition}

\section{Certified local-to-global reconstruction in \texorpdfstring{$W^{1,2}(0,1)$}{W12(0,1)}}\label{sec:reconstruction}

The remaining sections specialize the general evidence layer to an exact finite arithmetic
model.  The analytic estimates are standard; the additional object is the compiler that
turns supplied finite witnesses into a checkable global derivation with explicit resource
bounds.

\subsection{The rational piecewise-polynomial presentation}

Throughout Sections~\ref{sec:reconstruction}--\ref{sec:descent} we use the Hilbert norm convention
\[
 \norm{u}_{W^{1,2}(I)}^2
 :=\norm{u}_{L^2(I)}^2+\norm{u'}_{L^2(I)}^2,
 \qquad
 \langle u,v\rangle_{W^{1,2}(I)}
 :=\int_I uv+u'v'.
\]
All componentwise estimates below are combined relative to this convention.

Let $I=(a,b)$ have rational endpoints.

\begin{definition}[Rational $W^{1,2}$ presentation]\label{def:W12}
The finite code set $\Code_I$ consists of a rational mesh
\[
 a=t_0<t_1<\cdots<t_N=b
\]
and rational-coefficient polynomials $P_k\in\Q[x]$ on each
$[t_{k-1},t_k]$, with continuity at the interior mesh points.  Continuity is an exact
rational identity $P_k(t_k)=P_{k+1}(t_k)$, so $\Code_I$ is decidable.  The decoder sends
the code to the corresponding continuous piecewise-polynomial function in $W^{1,2}(I)$;
these codes are dense in $W^{1,2}(I)$.

A represented name $\nu$ is a fast Cauchy procedure producing codes $\nu(n)$ and, on
demand for each pair $m\ge n$, finite verification data for
\[
 \norm{\rho_I(\nu(m))-\rho_I(\nu(n))}_{W^{1,2}(I)}\le2^{-n}.
\]
Hence
\[
 \norm{\delta_I(\nu)-\rho_I(\nu(n))}_{W^{1,2}(I)}\le2^{-n}.
\]

For approximation evidence, the exact witness is accepted at bound $0$ for a canonical
constant name.  A positive witness $(n,\sigma)$ at rational bound $q$ is accepted when
$q>2^{-n}$,
\[
 \sigma=\norm{\rho_I(\nu(n))-\rho_I(p)}_{W^{1,2}(I)}^2\in\Q,
\]
and
\[
 \sigma<(q-2^{-n})^2.
\]
Distance evidence between two names is analogous: a witness at stage $n$ is accepted only
when $q>2^{1-n}$ and
\[
 \norm{\rho_I(\nu(n))-\rho_I(\mu(n))}_{W^{1,2}}^2
 <(q-2^{1-n})^2.
\]
The concrete evidence signature contains the normalized structural rules of
Definition~\ref{def:enrichment}, together with the declared multiplier, finite-sum,
bounded-overlap aggregation, and exact rational-arithmetic tags used below.
\end{definition}

\begin{lemma}[Exact finite-code arithmetic]\label{lem:exactW12}
For $p,q\in\Code_I$,
\[
 \langle\rho_I(p),\rho_I(q)\rangle_{W^{1,2}(I)}
 \quad\text{and}\quad
 \norm{\rho_I(p)-\rho_I(q)}_{W^{1,2}(I)}^2
\]
are rational numbers computable exactly by finitely many rational operations.  The code
language is effectively closed under finite sums, rational scaling, restriction to
rational subintervals, and multiplication by continuous rational piecewise-polynomial
codes.
\end{lemma}

\begin{proof}
Refine the two meshes to the union of their rational breakpoints; if the inputs have
$N_p,N_q$ cells, the common refinement has at most $N_p+N_q$ cells.  On each cell both functions and
their weak derivatives are rational polynomials.  Products have rational coefficients,
and their integrals over rational endpoints are rational.  Summing the cell integrals gives exact $L^2$ and derivative contributions.  For fixed
polynomial degree, addition and multiplication increase coefficient bit length by at most
a constant multiple of the sum of the input bit lengths (plus logarithmic finite-sum
overhead).  The closure statements follow from the same common-mesh operations and are
the exact finite arithmetic used by the synthesis compiler of Section~\ref{sec:onecover}.
\end{proof}

\begin{proposition}[Concrete soundness and strict-slack completeness]\label{prop:W12sound}
The presentation of Definition~\ref{def:W12} is sound and strict-slack complete on
$D_I=W^{1,2}(I)$.  Every finite code has exact principal evidence, and the evidence
regularity required in Proposition~\ref{thm:nonfaithful} holds.
\end{proposition}

\begin{proof}
For approximation evidence, acceptance gives
\[
 \norm{\rho_I(\nu(n))-\rho_I(p)}<q-2^{-n};
\]
the Cauchy tail gives the required true distance $<q$.  The distance checker is identical
with two tails.  Conversely, if the true distance is strictly below $q$, sufficiently large $n$ leaves positive rational slack, and Lemma~\ref{lem:exactW12}
computes the exact finite-stage squared distance.  Unlike Proposition~\ref{prop:slack}, no
rational-enclosure subsearch is needed here: the finite-stage squared norm is itself an
exact rational and the checker compares rationals directly.  Constant names give exact principal evidence.  Promotion
from an approximation witness to distance from the corresponding constant name uses the
same finite-stage arithmetic.
\end{proof}

\subsection{Rational partitions and synthesis}

\begin{lemma}[Certified rational partitions of unity]\label{lem:rationalPOU}
Let $\mathcal U=\{U_i\}_{i=1}^M$ be a finite cover of $I$ by relatively open rational
intervals.  One can effectively compute a rational subdivision fine enough that the closed star assigned to
each vertex lies inside some $U_i$, and from it compute a continuous rational piecewise-linear
partition of unity $(\psi_i)$ subordinate to $\mathcal U$ with
\[
 0\le\psi_i\le1,
 \qquad \sum_i\psi_i=1,
 \qquad
 L_i:=\norm{\psi_i'}_\infty\in\Q_{\ge0}.
\]
\end{lemma}

\begin{proof}
Because the cover is finite with rational interval endpoints, a positive rational mesh
scale subordinate to the cover can be found effectively by finite comparison of endpoints.
Use the usual rational nodal hat functions on a finer subdivision and assign each vertex to a patch containing
its closed star.  Summing the hats assigned to $U_i$ gives $\psi_i$.  Breakpoints,
heights, slopes, and derivative bounds are rational and effectively computable.
\end{proof}

For $v\in W^{1,2}(U_i)$, multiplication by $\psi_i$ followed by zero extension to $I$ is
a bounded linear map; the subordinate construction makes $\psi_i$ vanish at the relevant
boundary of its support, so the zero extension remains in $W^{1,2}(I)$.  Writing
$C_\infty=\norm{\psi_i}_\infty$ and $L_i=\norm{\psi_i'}_\infty$, the product rule gives
\[
 \norm{\psi_i v}_{W^{1,2}}^2
 \le (C_\infty^2+2L_i^2)\norm v_{L^2}^2
     +2C_\infty^2\norm{v'}_{L^2}^2.
\]
Hence one may choose any rational $C_i$ strictly above
$\sqrt{\max\{C_\infty^2+2L_i^2,2C_\infty^2\}}$.  On finite rational codes the action
is exact by Lemma~\ref{lem:exactW12}.  On a quasi-uniform hat mesh of scale $h$,
$L_i=O(h^{-1})$; consequently the local $L^2$ error $O(h^r)$ is amplified to
$O(h^{r-1})$, matching the derivative error and explaining the $r\ge2$ refinement regime.
The local multiplication/zero-extension fact and the resulting PUFEM estimates are
standard~\cite{MelenkBabuska96,BabuskaMelenk}.

\begin{proposition}[Exact finite synthesis]\label{prop:exactsynth}
Given local finite codes $p_i\in\Code_{U_i}$ and a rational partition $(\psi_i)$, the
finite global code
\[
 p_G:=\sum_i \psi_i p_i
\]
can be computed exactly in $\Code_I$.  Given local evidence for
$\norm{x_i-\rho(p_i)}_{W^{1,2}(U_i)}\le\alpha_i$, the declared multiplier and finite-sum
rules compile these witnesses into a global approximation proof for
$G=\sum_i\psi_i x_i$; for example the elementary triangle estimate gives
\[
 \norm{G-\rho(p_G)}_{W^{1,2}(I)}\le\sum_i C_i\alpha_i.
\]
This proposition synthesizes supplied local sections; relating $G$ to a common analytic
target is the separate PUFEM step of Theorem~\ref{thm:pufemdefect}.
\end{proposition}

\begin{proof}
Exactness of $p_G$ follows from the multiplication and finite-sum closure in
Lemma~\ref{lem:exactW12}.  Apply the stored multiplier rule to each local witness, extend
by zero, and combine the resulting accepted bounds by the finite-sum/triangle rule.
\end{proof}

\subsection{Localized PUFEM defect bound}

Let the cover have pointwise overlap multiplicity at most $\kappa$, and suppose
\[
 \norm{\psi_i}_\infty\le C_\infty,
 \qquad
 \norm{\psi_i'}_\infty\le L_i.
\]
For neighboring patches write $i\sim j$ when their relevant supports overlap.

\begin{theorem}[Certified localized PUFEM defect]\label{thm:pufemdefect}
Let $x_i\in W^{1,2}(U_i)$ be represented local sections with supplied finite evidence for
\[
 \norm{x_i-x_j}_{L^2(U_i\cap U_j)}\le\delta^0_{ij},
 \qquad
 \norm{(x_i-x_j)'}_{L^2(U_i\cap U_j)}\le\delta^1_{ij}.
\]
Define
\[
 G:=\sum_i\psi_i x_i
\]
with multiplication and zero extension understood patchwise.  For every $j$, finite
structural evidence can be compiled for
\[
 \norm{G|_{U_j}-x_j}_{W^{1,2}(U_j)}^2\le R_j,
\]
where
\[
 R_j
 =\kappa\sum_{i\sim j}
 \left[(C_\infty^2+2L_i^2)(\delta^0_{ij})^2
       +2C_\infty^2(\delta^1_{ij})^2\right].
\]
Equivalently, if
\[
 A_j^2=\kappa C_\infty^2\sum_{i\sim j}(\delta^0_{ij})^2,
\]
\[
 B_j^2=2\kappa\sum_{i\sim j}
 \left[L_i^2(\delta^0_{ij})^2+C_\infty^2(\delta^1_{ij})^2\right],
\]
then the componentwise squared bounds are
\[
 \norm{G|_{U_j}-x_j}_{L^2}^2\le A_j^2,
 \qquad
 \norm{(G|_{U_j}-x_j)'}_{L^2}^2\le B_j^2.
\]
For the rational partitions of Lemma~\ref{lem:rationalPOU} one may take $C_\infty=1$;
we retain $C_\infty$ to display the dependence for a general subordinate partition.
\end{theorem}

\begin{proof}
On $U_j$, the partition identity gives
\[
 G-x_j=\sum_{i\sim j}\psi_i(x_i-x_j).
\]
At every point at most $\kappa$ summands are nonzero, hence
\[
 \norm{G-x_j}_{L^2}^2
 \le\kappa\sum_{i\sim j}
       \norm{\psi_i(x_i-x_j)}_{L^2}^2
 \le A_j^2.
\]
For the weak derivative,
\[
 (G-x_j)'=\sum_{i\sim j}
 \left[\psi_i'(x_i-x_j)+\psi_i(x_i-x_j)'\right].
\]
Using bounded overlap and $(a+b)^2\le2a^2+2b^2$ gives
\[
 \norm{(G-x_j)'}_{L^2}^2
 \le2\kappa\sum_{i\sim j}
 \left[L_i^2(\delta^0_{ij})^2
 +C_\infty^2(\delta^1_{ij})^2\right]=B_j^2.
\]
Adding the squared $L^2$ and derivative estimates yields $R_j=A_j^2+B_j^2$.
The proof DAG is obtained by tagging the supplied local defect witnesses and applying the
declared multiplier, bounded-overlap, finite-sum, and rational-arithmetic rules.  If
$I=\sum_j\deg(j)$ is the number of oriented neighbour incidences, this construction uses
$O(I+M)$ new bounded-arity nodes before bit lengths are counted.  The
analytic inequality is the localized standard PUFEM argument; the additional assertion is
that the finite incoming evidence is sufficient to build the declared target derivation.
\end{proof}

\begin{remark}[Refinement-aware defect]
If $L_i\le C_Gh_i^{-1}$, the natural scale-sensitive local defect is
\[
 e_{ij}^2:=h_i^{-2}(\delta^0_{ij})^2+(\delta^1_{ij})^2.
\]
Only the zeroth-order mismatch is amplified by inverse patch scale.  Standard local approximation estimates supply one additional power of $h$ to the $L^2$ term:
$L_i=O(h_i^{-1})$ multiplying an $O(h_i^r)$ local $L^2$ error gives
$O(h_i^{r-1})$, the same scale as the derivative error.  This is the scale mechanism behind
the standing assumption $r\ge2$ in the refinement theorem.
\end{remark}

\section{Quantitative proof transport under one-cover gluing}\label{sec:onecover}

\subsection{The fixed proof-DAG encoding}

The following cost statements are meaningful only after an encoding is fixed.  The key
choice is a self-contained shared proof DAG: each supplied source certificate is stored once,
repeated subderivations are referenced rather than copied, and later nodes refer to earlier
nodes by internal binary indices.  Static analytic facts (triangle,
stored multiplier inequalities, bounded-overlap inequalities, and the PUFEM estimate)
are identified by fixed theorem tags.  Rational quantities carry their ordinary binary
numerator/denominator encoding, and exact piecewise-polynomial coefficients remain
ordinary rationals.  No representation-invariance claim is attached to this cost model.

\begin{definition}[Incremental transport resources]\label{def:incremental}
For an input evidence DAG $B_{\rm in}$ and an output evidence DAG $B_{\rm out}$, let
\[
 S(B)=\text{encoded bit length},
 \qquad
 \Delta S=S(B_{\rm out})-S(B_{\rm in}).
\]
Let $V$ denote checker bit-time, $T$ the compiler bit-time, and $Q$ the maximum represented
name/oracle lookahead.  Certificate generation from the underlying semantic point is not
counted as transport; source certificates are inputs to the theorem.
\end{definition}

Let
\[
 I:=\sum_j\deg(j)
\]
be the number of oriented overlap incidences.  Let $N_{\rm in}$ be the input DAG node count,
let $B$ bound the bit length of the quantitative rational input data, and let $s$ be a
requested dyadic slack precision for final square-root/enclosure operations.

\begin{theorem}[Single-cover provenance compiler]\label{thm:singlecovercost}
Under the hypotheses of Theorem~\ref{thm:pufemdefect}, the localized global-defect proofs
can be compiled evidence-locally with
\[
 \Delta S_{\rm defect}
 =O\!\left((I+M)\bigl[B+s+\log(N_{\rm in}+I+M+2)\bigr]\right).
\]
If $\mathcal A(b)$ bounds the bit-time for the fixed-arity exact rational operations on
$b$-bit arguments (for example, $\mathcal A(b)=O(b^2)$ is a coarse schoolbook bound), then
\[
 V_{\rm out}
 \le V_{\rm in}
 +O\!\left((I+M)\,\mathcal A(B+s+\log(I+M+2))\right).
\]
The compiler makes no independent query to a semantic name of the synthesized target:
\[
 Q_{\rm target}=0.
\]
For bounded scalar precision in a word model, the structural proof-node overhead is
$O(I+M)$ words.
\end{theorem}

\begin{proof}
For each oriented incidence $(i,j)$, the compiler references the supplied $L^2$ and
derivative defect witnesses once and appends a constant number of multiplier and product-
rule nodes.  For every target patch $j$, it appends a bounded-overlap aggregation node,
computes the rational $R_j$ of Theorem~\ref{thm:pufemdefect}, and, if a norm rather than a
squared norm is requested, computes a dyadic rational $q_j$ with
\[
 \sqrt{R_j}<q_j\le\sqrt{R_j}+2^{-s}.
\]
No search over target names occurs.  There are $O(I+M)$ new nodes.  A node label stores
$O(B+s)$ quantitative bits and $O(\log(N_{\rm in}+I+M+2))$ bits for its bounded number of
references.  Summation gives the size bound.  Verification traverses the new nodes once and performs $O(I+M)$ exact fixed-arity
rational operations, giving the stated time bound.
\end{proof}

General semantics for proof evidence and producer/checker architectures are established in
formal verification~\cite{ChihaniMillerRenaud,Dandelion}.  The role of the present compiler
is to connect that proof representation to explicit PUFEM analytic inputs, output bounds,
and resource accounting.

\subsection{Weighted precision allocation for exact synthesis}

Suppose local source certificates provide finite codes $p_i$ with
\[
 \norm{x_i-\rho(p_i)}_{W^{1,2}(U_i)}<\alpha_i.
\]
Because synthesis is exact on codes, put
\[
 p_G=\sum_i\psi_i p_i.
\]

\begin{proposition}[Weighted global approximation budget]\label{prop:weightedbudget}
With $G=\sum_i\psi_i x_i$,
\[
 \norm{G-\rho(p_G)}_{W^{1,2}}^2
 <\kappa\sum_i w_i^2\alpha_i^2,
 \qquad
 w_i^2:=\max\bigl\{C_\infty^2+2L_i^2,\ 2C_\infty^2\bigr\}.
\]
Hence every rational allocation satisfying
\[
 \kappa\sum_iw_i^2\alpha_i^2<\varepsilon^2
\]
produces a valid global $\varepsilon$-certificate.
\end{proposition}

\begin{proof}
Put $e_i=x_i-\rho(p_i)$ and abbreviate
$A_i:=\norm{e_i}_{L^2}^2$ and $D_i:=\norm{e_i'}_{L^2}^2$.  Then
$G-\rho(p_G)=\sum_i\psi_i e_i$, and bounded overlap gives
\[
 \norm{\sum_i\psi_i e_i}_{L^2}^2
 \le\kappa C_\infty^2\sum_i A_i.
\]
For the derivative, bounded overlap and $(a+b)^2\le2a^2+2b^2$ give
\[
 \norm{\sum_i(\psi_i'e_i+\psi_ie_i')}_{L^2}^2
 \le2\kappa\sum_i
 \left(L_i^2 A_i+C_\infty^2D_i\right).
\]
Adding the two components according to the fixed $W^{1,2}$ norm convention yields
\[
 \norm{G-\rho(p_G)}_{W^{1,2}}^2
 \le\kappa\sum_i\left[(C_\infty^2+2L_i^2)A_i+2C_\infty^2D_i\right].
\]
Since $A_i+D_i=\norm{e_i}_{W^{1,2}}^2<\alpha_i^2$, each summand is bounded by
\[
 \max\bigl\{C_\infty^2+2L_i^2,\ 2C_\infty^2\bigr\}\alpha_i^2,
\]
which proves the claim.
\end{proof}

If $C_{x_i}(\alpha)$ denotes the size of the selected local certificate at tolerance
$\alpha$, Proposition~\ref{prop:weightedbudget} yields the concrete carrying profile
\[
 C_G^{\rm carry}(\varepsilon)
 \le
 \inf_{\kappa\sum_iw_i^2\alpha_i^2<\varepsilon^2}
 \left[
   \sum_i C_{x_i}(\alpha_i)
   +S_{\rm syn}(p_1,\ldots,p_M;\psi)
   +O(M(B+\log(M+2)))
 \right].
\]
The infimum records a carrying profile rather than an algorithm: an actual compiler fixes
a rational feasible allocation.  The use of an error/work allocation under a global constraint is standard in numerical
analysis; the paper does not claim that optimization paradigm.  The additional term here
is the finite proof-carrying synthesis trace.

\begin{proposition}[Global rational-code size]\label{prop:codesize}
Let local code $p_i$ have $m_i$ mesh cells and polynomial degree at most $d_i$, let
$\psi_i$ have $r_i$ cells, and let $B$ bound input rational coefficient/endpoint bit
lengths.  If at most $\kappa$ partition terms are nonzero on any global cell, then one may
encode $p_G$ with
\[
 m_G\le\sum_i(m_i+r_i),
 \qquad
 d_G\le d_{\max}+1,
\]
and with a safe coefficient-bit bound
\[
 b_G=O\!\left(\kappa[B+\log(\kappa+1)]\right).
\]
Thus for bounded degree and bounded overlap the synthesized finite-code size is linear in
the total local/partition mesh volume, up to rational bit-length factors.
\end{proposition}

\begin{proof}
Multiplication by a piecewise-linear $\psi_i$ refines to the union of the two meshes and
raises polynomial degree by at most one.  A common global mesh therefore has at most the
sum of the local cell counts.  On each cell at most $\kappa$ exact rational polynomial
coefficients are added.  Standard numerator/denominator arithmetic gives the stated safe bit-length bound.  In the
refinement theorem each level is synthesized from fresh level-local codes, not recursively
from the previous global code, so the fixed overlap factor $\kappa$ is incurred once per
level and never compounds across refinement depth.
\end{proof}

\section{Proof-carrying refinement and represented-limit transport}\label{sec:descent}

We now turn the one-cover compiler into a refinement theorem.  The goal is not to reprove
constructive completeness.  It is to retain a finite ancestry of the precision-relevant
refinement process while controlling the additional proof resources.

\subsection{Provenance-complete certificates}

\begin{definition}[Provenance-complete certificate]\label{def:pcert}
A provenance-complete certificate is
\[
 \PCert=(p,\bar\varepsilon,V,H),
\]
where $(p,\bar\varepsilon,V)$ is an ordinary accepted certificate and $H$ is a finite
rooted DAG.  Every leaf of $H$ is an accepted input certificate or certified structural
datum; every internal node is an application of a declared finite evidence constructor;
and, with edges oriented from premises to conclusions, the distinguished root is the sink identifying the construction of $p$.  Shared subderivations are
stored once and referenced by binary index.  Its encoded size, using the DAG encoding of Section~\ref{sec:onecover}, is
\[
 |\PCert|=|p|+\ell_{\Q}(\bar\varepsilon)+|V|+|H|.
\]
The provenance checker accepts exactly when the ordinary certificate and every reachable
DAG node check.
\end{definition}

This enrichment changes no analytic truth condition.  It makes a complete finite
construction genealogy a checkable object rather than informal metadata.

\subsection{Scale-sensitive analytic input}

\begin{theorem}[Classical scale-sensitive PUFEM estimate]\label{thm:scalePUFEM}
Let $\mathcal V_h=\{V_a\}$ be a quasi-uniform interval cover of $(0,1)$ with patch diameter
comparable to $h\le1$, overlap at most $\kappa$, and a subordinate partition $(\chi_a)$
with
\[
 \norm{\chi_a'}_\infty\le C_\chi h^{-1}.
\]
Let $f\in W^{r,2}(0,1)$ with $r\ge2$, and suppose local approximants $u_a$ satisfy
\[
 \norm{f-u_a}_{L^2(V_a)}
 \le C_0h^r|f|_{W^{r,2}(V_a)},
\]
\[
 \norm{f'-u_a'}_{L^2(V_a)}
 \le C_1h^{r-1}|f|_{W^{r,2}(V_a)}.
\]
Then for $u_h=\sum_a\chi_au_a$,
\[
 \norm{f-u_h}_{W^{1,2}(0,1)}
 \le C_*h^{r-1}|f|_{W^{r,2}(0,1)},
\]
where one may take
\[
 C_*:=\kappa\bigl((1+C_\chi)C_0+C_1\bigr).
\]
\end{theorem}

\begin{proof}
Bounded overlap gives
\[
 \norm{f-u_h}_{L^2}
 \le\sqrt\kappa\left(\sum_a\norm{f-u_a}_{L^2(V_a)}^2\right)^{1/2}.
\]
Because $\sum_a\chi_a'=0$,
\[
 f'-u_h'=\sum_a\bigl[\chi_a'(f-u_a)+\chi_a(f'-u_a')\bigr].
\]
A second bounded-overlap estimate gives
\[
 \norm{f'-u_h'}_{L^2}
 \le\sqrt\kappa\left[
 C_\chi h^{-1}\left(\sum_a\norm{f-u_a}_{L^2}^2\right)^{1/2}
 +\left(\sum_a\norm{f'-u_a'}_{L^2}^2\right)^{1/2}
 \right].
\]
Finally $\sum_a|f|_{W^{r,2}(V_a)}^2\le\kappa|f|_{W^{r,2}}^2$.  Therefore
\[
 \norm{f-u_h}_{W^{1,2}}
 \le \kappa h^{r-1}|f|_{W^{r,2}}
 \sqrt{C_0^2h^2+(C_\chi C_0+C_1)^2}
 \le \kappa((1+C_\chi)C_0+C_1)h^{r-1}|f|_{W^{r,2}},
\]
using $h\le1$ and $\sqrt{a^2+b^2}\le a+b$.  This is standard partition-of-unity
approximation theory~\cite{MelenkBabuska96,BabuskaMelenk}.
\end{proof}

\subsection{Same-order genealogy transport}

\begin{proposition}[Concrete compiler bound for standard local encodings]\label{prop:H6compiler}
Consider a geometric refinement hierarchy with $M_n$ level-$n$ patches and
$M_j\le C_g2^{j-n}M_n$ for $j\le n$, together with uniformly bounded same-level
neighbour degree and refinement incidence.  Suppose each local approximation and partition
code has $O(1)$ mesh cells and uniformly bounded polynomial degree; all endpoints,
coefficients, supplied local witnesses, and proof-node labels use $O(\beta_n)$ bits;
the root-slack schedule satisfies $s_n=O(\beta_n)$; and
$\beta_n\ge c\log_2(M_n+2)$.  Assume also that the conventional non-proof baseline is
nondegenerate,
\[
 B_n=\Theta(M_n\beta_n).
\]
Then the concrete compilers of Sections~\ref{sec:reconstruction}--\ref{sec:onecover}
produce
\[
 |p_n|=O(M_n\beta_n)=O(B_n)
\]
and introduce only $O(M_n)$ new bounded-arity nodes at level $n$, with total new proof
payload $O(M_n\beta_n)$.  Consequently the accumulated node count before level $n$ is
also $O(M_n)$.
\end{proposition}

\begin{proof}
The bounded neighbour hypothesis gives $I_n=O(M_n)$.  Proposition~\ref{prop:codesize}
gives $O(M_n)$ global cells, bounded degree, and coefficient bit length
$O(\kappa[\beta_n+\log(\kappa+1)])=O(\beta_n)$ because $\kappa$ is a fixed geometry
constant.  Hence $|p_n|=O(M_n\beta_n)$.  The one-cover construction introduces
$O(I_n+M_n)=O(M_n)$ new bounded-arity nodes.  Geometric growth therefore gives an input
DAG node count
\[
 N_{{\rm in},n}=O\!\left(\sum_{j<n}M_j\right)=O(M_n).
\]
Theorem~\ref{thm:singlecovercost} now yields
\[
 \Delta S_n
 =O\!\left(M_n[\beta_n+s_n+\log(M_n+2)]\right)
 =O(M_n\beta_n),
\]
because its reference term
$\log(N_{{\rm in},n}+I_n+M_n+2)$ is $O(\log(M_n+2))=O(\beta_n)$.
The supplied local witnesses themselves contribute $O(M_n\beta_n)$ by assumption, and
the bounded refinement incidences add only $O(M_n)$ further bounded-arity references and
labels of the same bit scale.  This proves the stated node and payload bounds.
\end{proof}

The proposition deliberately separates the conventional baseline from the proof-carrying
conclusion: $B_n=\Theta(M_n\beta_n)$ and
$\beta_n\gtrsim\log(M_n+2)$ declare the encoding against which overhead is measured;
$|p_n|=O(B_n)$ and the $O(M_n\beta_n)$ proof payload are derived from the compiler.
A linear schedule $\beta_n=O(n+1)$ is a further property of the chosen local encoding,
not a consequence of Sobolev approximation alone.

\begin{theorem}[Encoding-cost evidence transport under Sobolev refinement]\label{thm:orderneutral}
Let $\Omega=(0,1)$, let $r\in\N$ with $r\ge2$, and put $\alpha=r-1$.
Fix $h_0\in\Q\cap(0,1]$ and $h_n=h_0 2^{-n}$.  For every $n\ge0$ let
\[
 \mathcal V_n=\{V_{n,a}:a\in A_n\},\qquad M_n:=|A_n|,
\]
be a rational interval cover.  Assume:
\begin{enumerate}[label=(H\arabic*)]
\item \textbf{Geometric quasi-uniformity.} There are constants $c_M,C_M>0$ and
integers $\kappa,D,D_R\ge1$, independent of $n$, such that
\[
 c_Mh_n^{-1}\le M_n\le C_Mh_n^{-1},
\]
pointwise overlap is at most $\kappa$, every patch has at most $D$ relevant same-level
neighbors, and the number of refinement incidences between levels $n$ and $n+1$ is at
most $D_RM_{n+1}$.

\item \textbf{Certified rational partitions.} Each cover carries a rational continuous
piecewise-linear partition $(\chi_{n,a})_{a\in A_n}$ subordinate to $\mathcal V_n$ with
\[
 0\le\chi_{n,a}\le1,\qquad
 \sum_a\chi_{n,a}=1,\qquad
 \norm{\chi_{n,a}'}_\infty\le C_\chi h_n^{-1},
\]
for one rational $C_\chi$ independent of $n$.

\item \textbf{Certified local order-$r$ approximation.} There is a represented
$f\in W^{r,2}(0,1)$ and a rational $R\ge |f|_{W^{r,2}}$.  For every $n,a$ there is a
rational piecewise-polynomial code $u_{n,a}$ and supplied accepted finite evidence for
\[
 \norm{f-u_{n,a}}_{L^2(V_{n,a})}
 \le C_0h_n^r |f|_{W^{r,2}(V_{n,a})},
\]
\[
 \norm{f'-u_{n,a}'}_{L^2(V_{n,a})}
 \le C_1h_n^{r-1}|f|_{W^{r,2}(V_{n,a})},
\]
where $C_0,C_1\in\Q_{\ge0}$ are independent of $n,a$.

\item \textbf{Exact synthesis and transition evidence.} The global finite code
\[
 p_n:=\sum_{a\in A_n}\chi_{n,a}u_{n,a}
\]
is compiled exactly.  The level-$n$ derivation includes accepted finite evidence for the
scale estimate
\[
 \norm{f-\rho(p_n)}_{W^{1,2}}\le E_n,
 \qquad
 E_n:=C_*R h_n^{\alpha},
\]
with $C_*$ as in Theorem~\ref{thm:scalePUFEM}.  Adjacent transition evidence is obtained
from the two comparison witnesses by symmetry and triangle.

\item \textbf{Persistent evidence-local genealogy.} A uniform compiler constructs finite
DAGs $H_n$ such that $H_{n+1}$ contains $H_n$ by shared reference, not by copying.  The
new level records the local input witnesses, finite synthesis, scale estimate, and the
transition from level $n$ to $n+1$.  The compiler never invokes an independent generic
certificate generator on the eventual semantic limit.

\item \textbf{Declared coefficient-wise encoding.} There is a nondecreasing sequence
$\beta_n\ge1$ and positive constants $b_-,b_+,b_p,b_H,b_N$ such that
\[
 \beta_n\ge b_H^{-1}\log_2(M_n+2),
\]
and the conventional non-proof encoding size $B_n$ of the level-$n$ finite approximation
data satisfies
\[
 b_-M_n\beta_n\le B_n\le b_+M_n\beta_n,
 \qquad
 |p_n|\le b_pB_n.
\]
The number of new bounded-arity proof nodes introduced at level $n$ is at most $b_NM_n$,
and the total bit length of all \emph{new} local evidence payloads and proof-node labels is
at most $b_HM_n\beta_n$.  DAG references are binary indices.  Proposition~\ref{prop:H6compiler}
shows that these proof-carrying clauses follow from transparent standard local-encoding
assumptions; the theorem keeps (H6) abstract so other encodings can be used.

\item \textbf{Arithmetic verification model.} Every fixed-arity rational operation used
by a level-$n$ proof node acts on $O(\beta_n)$-bit quantities.  Let $\mathcal A(b)$ be a
nondecreasing upper bound for the bit-time required to check one such node on $O(b)$-bit
data.
\end{enumerate}
Then there is a computable name $\nu_\infty$ for $f$ and a provenance-complete certificate
system $c_\infty^{\rm pc}$ over that name such that, for rational
$0<\varepsilon\le1$,
\[
 |c_\infty^{\rm pc}(\varepsilon)|
 \le C B_{m(\varepsilon)}
\]
for one constant $C$ depending only on the fixed geometry, analytic, and encoding
parameters in (H1)--(H7), but not on $\varepsilon$.  The proof constructs a computable
level selector $m(\varepsilon)$ satisfying
\[
 m(\varepsilon)
 =O\!\left(\frac{\log(1/\varepsilon)}{r-1}\right).
\]
Moreover, the persistent history satisfies
\[
 \sum_{j\le n}M_j\beta_j=O(M_n\beta_n),
\]
so the complete genealogy through the precision-relevant level has the same asymptotic bit
order as the finest-level declared coefficient-wise encoding.

If, in addition, $\beta_n=O(n+1)$, then
\[
 \boxed{
 |c_\infty^{\rm pc}(\varepsilon)|
 =O\!\left(
 \varepsilon^{-1/(r-1)}
 (1+\log(1/\varepsilon))
 \right).}
\]
The construction uses no independently supplied semantic target name,
\[
 \boxed{Q_{\rm target}=0.}
\]
If, in addition, $\beta_n=O(n+1)$ and generation of the supplied level-$n$ local
certificates uses represented-name precision lookahead $O(\beta_n)$, then the additional
source-precision requirement of the refinement-to-limit construction satisfies
\[
 \boxed{Q_{\rm source}(\varepsilon)=O(\log(1/\varepsilon)).}
\]
Finally, apart from the already declared cost of checking primitive local witnesses, the
complete lineage can be verified in time
\[
 \boxed{
 V(\varepsilon)
 =O\!\left(
 M_{m(\varepsilon)}\,
 \mathcal A(\beta_{m(\varepsilon)})
 \right).}
\]
\end{theorem}

\begin{proof}
By Theorem~\ref{thm:scalePUFEM}, Hypotheses (H2)--(H4) give
\[
 E_n=C_*Rh_0^\alpha 2^{-\alpha n}=:E_0\rho^n,
 \qquad \rho:=2^{-\alpha}\in(0,1).
\]
The level-$n$ and level-$(n+1)$ comparison witnesses, followed by symmetry and triangle,
give finite transition evidence
\[
 \norm{\rho(p_{n+1})-\rho(p_n)}
 \le E_n+E_{n+1}
 =(1+\rho)E_n.
\]
Hence the certified remaining transition tail from level $n$ is bounded by
\[
 \Theta_n
 :=\sum_{j=n}^{\infty}(E_j+E_{j+1})
 =\frac{1+\rho}{1-\rho}E_n
 =:K_rE_n.
\]
This provides a computable Cauchy modulus directly from the stored ancestry.

Choose an integer $c_\Theta\ge0$ with
$K_rE_0\le2^{c_\Theta}$.  Define
\[
 \mu(s)
 :=\left\lceil\frac{s+c_\Theta+3}{\alpha}\right\rceil,
 \qquad
 \nu_\infty(s):=p_{\mu(s)}.
\]
For $t\ge s$, monotonicity of $\mu$ and the stored transition chain give
\[
 \norm{\rho(\nu_\infty(t))-\rho(\nu_\infty(s))}
 \le\Theta_{\mu(s)}
 \le2^{c_\Theta}2^{-\alpha\mu(s)}
 \le2^{-s-3}<2^{-s}.
\]
The finite triangle chain extracted from the genealogy is the verification data required
by the fast Cauchy representation.  Thus $\nu_\infty$ is a computable represented name.
Since
\[
 \norm{f-\rho(p_{\mu(s)})}\le E_{\mu(s)}\to0,
\]
its represented limit is $f$.

Now fix rational $0<\varepsilon\le1$ and put
\[
 k_\varepsilon:=\left\lceil\log_2\frac2\varepsilon\right\rceil,
 \qquad
 \bar\varepsilon:=2^{-k_\varepsilon},
 \qquad
 s_\varepsilon:=k_\varepsilon+2,
 \qquad
 n_\varepsilon:=\mu(s_\varepsilon).
\]
Then $0<\bar\varepsilon\le\varepsilon/2<\varepsilon$ and
$2^{-s_\varepsilon}=\bar\varepsilon/4<\bar\varepsilon$.  At name stage
$s_\varepsilon$ the code is exactly $p_{n_\varepsilon}$.  In the checker of
Definition~\ref{def:W12}, the finite-stage squared distance between
$\nu_\infty(s_\varepsilon)$ and $p_{n_\varepsilon}$ is zero.  Hence the positive
witness $(s_\varepsilon,0)$ is accepted at the dyadic announced bound
$\bar\varepsilon$.  Attach the genealogy
$H_{n_\varepsilon}$ to obtain the provenance-complete certificate.  Notice that this
certificate is produced from the constructed refinement sequence itself; no separate
semantic target recertifier is called.

It remains to bound $|H_n|$.  By bounded overlap, bounded refinement incidence, bounded
proof-node arity, and Hypothesis (H6), the number and payload size of nodes newly introduced
at level $j$ are bounded by $c_NM_j$ and $c_SM_j\beta_j$, respectively.  Quasi-uniform
geometric refinement gives, for $j\le n$,
\[
 M_j\le\frac{C_M}{c_M}2^{j-n}M_n.
\]
Therefore
\[
 \sum_{j=0}^n M_j
 \le 2\frac{C_M}{c_M}M_n.
\]
Since $\beta_j\le\beta_n$,
\[
 |H_n|
 \le C_2\sum_{j=0}^nM_j\beta_j
 \le C_3M_n\beta_n.
\]
The binary index length of a DAG reference is also $O(\log(M_n+2))=O(\beta_n)$ by
Hypothesis (H6), so it is already absorbed in this estimate.  Since $B_n\ge b_-M_n\beta_n$ and $|p_n|\le b_pB_n$, the code and genealogy
terms are $O(B_n)$.  The dyadic announced slack and ordinary checker witness in the final certificate have
size $O(\log(1/\varepsilon))$, which is dominated by
$B_{n_\varepsilon}=\Omega(M_{n_\varepsilon})$.  Hence
\[
 |c_\infty^{\rm pc}(\varepsilon)|
 \le C B_{n_\varepsilon}
\]
for a fixed $C$ with the dependence stated in the theorem.

Next,
\[
 n_\varepsilon
 \le\frac{s_\varepsilon+c_\Theta+3}{r-1}+1
 =O\!\left(\frac{\log(1/\varepsilon)}{r-1}\right).
\]
Thus we take $m(\varepsilon):=n_\varepsilon$ in the theorem statement.  Also $M_n\le C_Mh_0^{-1}2^n$, hence
\[
 M_{n_\varepsilon}
 =O\!\left(\varepsilon^{-1/(r-1)}\right).
\]
If $\beta_n=O(n+1)$, then
$\beta_{n_\varepsilon}=O(1+\log(1/\varepsilon))$, which yields the displayed explicit
certificate-size estimate.

For verification, the total number of nodes through level $n_\varepsilon$ is
$O(M_{n_\varepsilon})$, and all labels have bit scale at most
$O(\beta_{n_\varepsilon})$.  Monotonicity of $\mathcal A$ therefore gives the stated
verification-time bound.  The target-query statement is immediate from the explicit
construction of $\nu_\infty$.  Finally, if primitive local certificate generation at level
$n$ needs source-name precision $O(\beta_n)$, the deepest queried level is
$n_\varepsilon$, so the additional source-precision lookahead is
$O(\beta_{n_\varepsilon})=O(\log(1/\varepsilon))$ under the stated linear bit schedule.
\end{proof}

\begin{corollary}[Same-order transport in the standard rational regime]\label{cor:orderneutral}
Under the hypotheses of Theorem~\ref{thm:orderneutral} with bounded polynomial degree,
bounded overlap, and ordinary binary rational coefficients whose per-patch precision grows
linearly with refinement depth,
\[
 |c_\infty^{\rm pc}(\varepsilon)|
 =O\!\left(B_{m(\varepsilon)}\right)
 =O\!\left(\varepsilon^{-1/(r-1)}(1+\log(1/\varepsilon))\right).
\]
The first comparison is the same-order statement: the complete proof-carrying package is
bounded by a constant multiple of the declared conventional finite encoding at the
precision-relevant level.  Linear growth of rational precision with refinement depth is a standard regime: for
bounded coefficient magnitudes, absolute dyadic resolution at an $h_n^r$ scale uses only
$O(n)$ fractional bits when $h_n\asymp2^{-n}$.  The corollary assumes that the chosen
local approximation scheme realizes such a schedule; no general quantization theorem is
claimed here.  No two-sided lower bound for certificate size is asserted.
\end{corollary}

Classical entropy-number theory for Sobolev and related function-space embeddings already
captures the intrinsic algebraic exponent~\cite{EdmundsTriebel}.  Mao and Xu compare
bit-complexity rates for neural-network and classical approximation schemes, explicitly
including finite elements, and distinguish coefficient-wise encodings from entropy-oriented
descriptions~\cite{MaoXu2026}.  The theorem here addresses a different quantity: whether the
additional finite proof genealogy introduces an asymptotic overhead beyond the declared
ordinary encoding.

\subsection{Why the resource hypotheses are substantive}

The same-order conclusion can fail when the hypotheses that make historical ancestry
summable, or that compare proof payload with ordinary data, are removed.

\begin{example}[Failure without summable refinement history]\label{ex:nonsummable}
Suppose the finest-level ordinary encoding still satisfies $B_n=\Theta(M_n\beta_n)$, but
replace geometric mesh growth by $M_n=n+1$ and take
\[
 \beta_n:=\lceil\log_2(n+3)\rceil\,.
\]
Then the index-length requirement in (H6) is satisfied (after changing its fixed constant
if necessary).  If each level contributes $\Theta(M_n\beta_n)$ genuinely new persistent
proof payload, as allowed by the same-level part of (H6), then
\[
 B_n=\Theta(n\log n),\qquad
 |H_n|=\Theta\!\left(\sum_{j=0}^n(j+1)\log(j+3)\right)=\Theta(n^2\log n).
\]
Hence $|H_n|/B_n=\Theta(n)\to\infty$.  Thus this example keeps the payload
comparability of (H6) while breaking the geometric-growth part of (H1).  It isolates the failure of the geometric-summability
step $\sum_{j\le n}M_j=O(M_n)$ used in Theorem~\ref{thm:orderneutral}: the bit schedule and
same-level proof-payload comparison remain admissible, but historical ancestry is no longer
summable relative to the finest level.  An adaptive or otherwise nongeometric refinement
history therefore needs its own cumulative-incidence hypothesis.
\end{example}

\begin{example}[Failure without proof-payload comparability]\label{ex:payload}
Keep geometric refinement $M_n=\Theta(2^n)$ and an ordinary rational precision schedule
$\beta_n=\Theta(n+1)$, so $B_n=\Theta(2^n(n+1))$.  If, contrary to (H6), the new evidence
payload per patch at level $n$ has bit scale $\gamma_n=\Theta(2^n)$, then the newest level
alone contributes
\[
 \Theta(M_n\gamma_n)=\Theta(4^n),
\]
and hence the persistent genealogy cannot be $O(B_n)$.  This example keeps geometric growth (H1) and breaks the payload comparison in (H6).
That comparison is therefore not bookkeeping decoration: it is what ties the cost of retained proof evidence
to the cost scale against which same-order transport is asserted.
\end{example}

These examples also clarify what the theorem does and does not say.  It is not the abstract
fact that costs add under composition.  The $h^{r-1}$ scale law fixes $m(\varepsilon)$; bounded
overlap controls same-level dependencies; geometric mesh growth makes earlier ancestry
summable relative to the finest level; and (H6) prevents proof payload from outgrowing the
declared approximation data.  Resource monoids can model the resulting bound, but they do not derive these PUFEM-specific
conditions~\cite{HofmannLedent}.  The scale hypothesis is load-bearing as well: if $r=1$,
then $\alpha=0$ and $\rho=1$, so the certified transition tail used to construct the Cauchy
modulus does not converge.  These are hypothesis-necessity examples, not lower bounds on
the size of every possible certificate.

\section{Contextual Choice: generated genealogies and admissible global construction}\label{sec:ccp}

Contextual Choice is used here as an operational discipline for a declared certified
development.  Its decisive requirement is that convergence, selection, reconstruction,
and local-to-global promotion enter only through explicit moduli or other declared finite
verification interfaces; computability alone is not taken to supply those interfaces.

\begin{definition}[CCP-admissible constructor]\label{def:ccpadm}
Let $\mathcal X$ be a declared evidence input type: a finite product of evidence
categories, or a certified sequence type equipped with explicit moduli and finite-query
interfaces sufficient to produce each requested finite output certificate.  A map
\[
 \Phi:\mathcal X\longrightarrow\Cert_{\rm ev}(\mathfrak Y)
\]
is \emph{CCP-admissible} when:
\begin{enumerate}[label=(\roman*)]
\item there is an underlying analytic operation $\phi$ satisfying
$U_Y\Phi=\phi U_{\mathcal X}$;
\item at every rational output tolerance, the output certificate is produced by a
terminating uniform algorithm from finitely many supplied input certificates and finite
structural data;
\item the algorithm returns finite evidence accepted by the target checker;
\item every convergence, selection, or reconstruction hypothesis used by the algorithm is
supplied with an explicit computable modulus or another declared finite verification
interface.
\end{enumerate}
\end{definition}

The qualitative notion deliberately does not fix a resource class.  A \emph{quantitative
CCP-admissible} constructor may additionally be required to satisfy declared bounds in the
resource vector $(G,V,T,S,Q;L)$ of Definition~\ref{def:resources}.

\begin{proposition}[Closure under admissible constructors]\label{prop:ccpclosure}
CCP-admissible constructors are closed under identities, finite products, and composition.
Finite-code realizable Lipschitz maps are admissible by Theorem~\ref{thm:saturation}.
Certified finite reconstruction is admissible by Proposition~\ref{prop:exactsynth} and
Theorem~\ref{thm:pufemdefect}; refinement sequences satisfying
Theorem~\ref{thm:orderneutral} admit a CCP-admissible limit constructor.
\end{proposition}

\begin{proof}
Identities copy evidence.  Finite products use tuple coding and coordinatewise checking.
For composition, run the first finite evidence transformer and then the second; the output
proof DAG contains the first output by reference and appends the second derivation.
Termination, finite evidence, and the underlying analytic commuting equation are
preserved.  The remaining claims are the cited constructions.
\end{proof}

\begin{definition}[Contextual Choice Principle]\label{def:ccp}
A certified development \emph{obeys CCP} when every claimed construction is either a
declared primitive or the output of a stated CCP-admissible constructor, and every
local-to-global promotion states the compatibility, reconstruction, refinement, and
modulus hypotheses used to justify it.
\end{definition}

The word ``choice'' now refers to admissible promotion of contextual local selections.  If
local selectors $s_i$ are supplied, CCP does not infer a global selector by fiat.  It
requires a proved operation of the form
\[
\begin{gathered}
 \text{certified local selections}+\text{certified compatibility}\\
 +\text{admissible gluing/refinement}
 \ \Longrightarrow\ 
 \text{certified global selection}.
\end{gathered}
\]
Nothing here forbids an external classical choice that lies outside the declared certified
signature.

\subsection{The generated universe}

\begin{definition}[CCP-generated universe]\label{def:generated}
Let $S$ be a set of evidence types.  For each $s\in S$, let $\mathcal E_s$ be the declared
evidence objects, $U_s:\mathcal E_s\to D_s$ the underlying analytic map,
$\mathcal P_s\subseteq\mathcal E_s$ a family of certified primitives, and $\Sigma$ a typed
signature of CCP-admissible constructors.  A typed family $\mathcal X=(\mathcal X_s)$ is
$\Sigma$-closed if it contains all primitives and is closed under every constructor in
$\Sigma$ whenever the constructor's declared side conditions hold.  Define
\[
 \mathsf{CU}_\Sigma(\mathcal P)
 :=\bigcap\{\mathcal X:\mathcal P\subseteq\mathcal X,
                 \ \mathcal X\text{ is }\Sigma\text{-closed}\}.
\]
\end{definition}

Here ``universe'' means only this least typed closure inside the declared evidence types;
it is neither a Grothendieck universe nor a claim to exhaust ordinary mathematical
objects.

\begin{theorem}[Invariant preservation]\label{thm:preservation}
For each type $s$, let $R_s\subseteq D_s$ be a regularity class.  Suppose every primitive
has underlying object in $R_s$ and every constructor in $\Sigma$ preserves the
corresponding regularity classes whenever its side hypotheses hold.  Then
\[
 U_s\bigl(\mathsf{CU}_\Sigma(\mathcal P)_s\bigr)\subseteq R_s
 \qquad(s\in S).
\]
Consequently every class $B_s\subseteq D_s$ disjoint from $R_s$ is absent from the
generated certified universe.
\end{theorem}

\begin{proof}
Let $\mathcal E_s^R=U_s^{-1}(R_s)$.  By the primitive assumption,
$\mathcal P_s\subseteq\mathcal E_s^R$, and preservation by every constructor says exactly
that the typed family $\mathcal E^R$ is $\Sigma$-closed.  Leastness of
$\mathsf{CU}_\Sigma(\mathcal P)$ gives
$\mathsf{CU}_\Sigma(\mathcal P)\subseteq\mathcal E^R$.
\end{proof}

The theorem is ordinary structural induction over a generated class; its role is semantic
clarity, not a priority claim.

\subsection{A finite-measure boundary example}

Let $\mathcal M$ be the metric space of measurable subsets of $\R^3$ of finite
Lebesgue measure, modulo null symmetric difference, with metric
\[
 d([A],[B])=\lambda^3(A\triangle B).
\]
Finite unions of bounded rational polyhedra give a countable dense core.  This finite-measure metric space is complete; effective limits additionally require a
represented Cauchy modulus.  Finite unions, intersections, set differences, and rational
Euclidean isometries preserve finite measure and have exact actions on the rational-polyhedral
core, and the measure is
$1$-Lipschitz with respect to symmetric difference.

\begin{corollary}[Internal Banach--Tarski boundary]\label{cor:BT}
Let the measurable-set component of $\mathsf{CU}_\Sigma(\mathcal P)$ be generated from
finite rational-polyhedral primitives by finite unions, intersections, and set differences,
certified measure-preserving Euclidean isometries, and effective limits in $\mathcal M$.  Then no finite certified equidecomposition inside this generated universe can
duplicate a set of finite positive measure.  In particular, no internal finite
equidecomposition can witness a Banach--Tarski duplication there.
\end{corollary}

\begin{proof}
Every admitted primitive is measurable of finite measure and every declared constructor
preserves that class, so Theorem~\ref{thm:preservation} keeps the generated objects in the
finite-measure class.  If
$A=\bigsqcup_{i=1}^N A_i$ and the admitted isometries $g_i$ produced disjoint images whose
union were two disjoint copies of $A$, finite additivity and isometry invariance would give
\[
 2\lambda^3(A)
 =\sum_i\lambda^3(g_iA_i)
 =\sum_i\lambda^3(A_i)
 =\lambda^3(A),
\]
contradicting $0<\lambda^3(A)<\infty$.
\end{proof}

This is an internal generated-class statement.  Its proof uses measurability and finite
additivity, not scarcity of the admitted isometries.  The classical Banach--Tarski theorem
uses nonmeasurable pieces and is untouched~\cite{BanachTarski1924}.

\section{Information boundaries and exact extensional gluing}\label{sec:boundaries}

Two boundary results help prevent the evidence language from being interpreted too
broadly.

\subsection{A model-relative information obstruction}

\begin{proposition}[Fibre indistinguishability]\label{prop:fibre}
Let $X$ be a normed vector space, let $I:X\to\R^m$ be linear and noninjective, and suppose a
sound verifier for a claimed finite norm bound $\norm x\le B$ can inspect $x$ only through
$I(x)$ and finite auxiliary data that do not distinguish points in the same fibre.  Then no finite-bound witness can be accepted on any nonempty fibre.
\end{proposition}

\begin{proof}
Every nonempty fibre is $x_0+\ker I$ and therefore contains an affine line
$x_0+\R z$ for some $z\in\ker I\setminus\{0\}$.  Every point on that line has the same information value $I(x_0)$.  If one finite witness were
accepted on that fibre, soundness under the stated information semantics would force the
same finite bound for all $\lambda$.  But
$\norm{x_0+\lambda z}\to\infty$ as $|\lambda|\to\infty$.
\end{proof}

This proposition concerns finite-dimensional linear information $I$, not the Type-2 name
access used in Sections~\ref{sec:evidence}--\ref{sec:descent}.  It is a standard
information-theoretic phenomenon and records why
certificate completeness always depends on what the verifier is allowed to read; compare
information-based complexity and finite-query SCI formulations~\cite{TWW,ColbrookHansen,Sorg2026}.

\subsection{Exact sheaf semantics}

Let $(\mathcal C,J)$ be a small certified site with Grothendieck topology $J$, whose
morphisms include certified restrictions, and on which represented sections satisfy exact
compatible gluing.  Define the extensional presheaf
\[
 \mathcal E(U):=\{\text{represented certified sections over }U\}/\text{zero evidence distance}.
\]

\begin{proposition}[Extensional sheaf placement]\label{prop:sheaf}
If exact compatible local sections glue uniquely in the represented analytic semantics and
certified restriction realizes the presheaf maps, then $\mathcal E$ is an ordinary sheaf
on $\mathcal C$.  For a small site it is therefore an object of the Grothendieck topos
$\mathbf{Sh}(\mathcal C,J)$.
\end{proposition}

\begin{proof}
By Theorem~\ref{thm:collapse}, after the zero-distance quotient matching families are
ordinary compatible represented sections.  Existence and uniqueness of their analytic gluing are exactly the ordinary
sheaf axiom.  The topos statement is standard sheaf theory~\cite{Johnstone,FourmanScott}.
\end{proof}

No stronger conclusion follows automatically.  Quantitative approximate gluing is not an
enriched sheaf universal property merely because it carries a numerical defect, and the
raw proof-relevant evidence assignment is not claimed to be a stack.

The categorical architecture used above requires no probes--models adjunction.  Its pieces
are represented information and finite codes, explicit realization/synthesis maps, Type-2
name access, chosen evidence lifts and their Grothendieck construction, proved
reconstruction and refinement, and ordinary extensional sheaf semantics.  Generic
certificate existence cannot itself define the admissible universe because
Theorem~\ref{thm:collapse} makes it extensionally cheap under full name access; the CCP
discipline therefore regulates construction genealogy instead.

\section{Relation to established frameworks}\label{sec:related}

\subsection{Represented spaces and realizability}

The represented-space layer is standard TTE~\cite{Weihrauch,BrattkaHertling}.  The
canonical slack argument of Proposition~\ref{prop:slack} is therefore infrastructure, not
an alternative foundation of computable analysis.  Realizability and assembly semantics
also already organize underlying functions by the existence of trackers~\cite{vanOosten}.
Our chosen-lift category retains a particular evidence transformer as arrow data; the
non-faithfulness proposition records the information lost by forgetting that extra choice.  It
does not imply a new extensional class of computable functions.

The same caution applies to finite-information lower bounds.  Information-based
complexity and SCI frameworks make the queried information interface explicit, and a recent
preprint emphasizes the danger of unrestricted post-processing of finite transcripts
\cite{TWW,ColbrookHansen,Sorg2026}.  Proposition~\ref{prop:fibre} is an illustration of
that established principle.

\subsection{Certified numerics and proof resources}

Certificate-producing rigorous approximation is also established.  Br\'ehard, Mahboubi, and Pous
develop formally verified approximation structures with a posteriori validation and compositional
rigorous operations~\cite{BrehardMahboubiPous}.  Dandelion provides a fully verified checker for
certificates attached to externally generated polynomial approximations~\cite{Dandelion}.  Certified
reduced-basis methods make rigorous a posteriori error bounds and guaranteed fidelity central to
reduced-order approximation~\cite{HesthavenRozzaStamm}; reconstruction-based a posteriori
finite-element methods likewise obtain locally computable guaranteed bounds through conforming
reconstructions~\cite{ErnStephansenVohralik}.  Formal finite-element infrastructure is an active area;
recent Rocq work formalizes simplicial Lagrange finite elements as records containing both mathematical
data and proofs of validity~\cite{BoldoEtAl2026}.  The paper therefore does not claim novelty for ``finite numerical object + proof'',
external producer plus trusted checker, or a posteriori certification as patterns.  The delta addressed here is different: supplied evidence is reused through assembly and
refinement without semantic target recertification, and the complete retained genealogy is
shown, under the stated encoding hypotheses, to have the same asymptotic bit order as the
bare finest-level approximation.  The cited frameworks establish rigorous final-object
certification and, in some cases, compositional certified operations; the present theorem
is specifically about quantitative retention of incoming proof genealogy through PUFEM
refinement.

Resource-sensitive semantics is likewise established.  Dal Lago and Hofmann introduced resource
monoids as abstract bounds in a realizability-based quantitative semantic framework~\cite{DalLagoHofmann};
Hofmann and Ledent use a simplified resource-monoid model to majorize normal-form size rather than
only execution time~\cite{HofmannLedent}.  Our $(G,V,T,S,Q;L)$ vector is not proposed as a replacement.
Its purpose is to state the concrete resources of one analytic proof-transport compiler.

Proof mining is another neighbouring programme: proof interpretations extract explicit
moduli, rates, and other quantitative data from existing proofs~\cite{Kohlenbach}.  The
present paper starts instead from supplied finite evidence and studies its transport and
retention cost through a fixed analytic construction.  The phrase ``proof-carrying'' also
recalls proof-carrying code, where an untrusted producer supplies a proof that a consumer
can check cheaply~\cite{NeculaPCC}; here the carried object is an analytic approximation
with a construction genealogy rather than executable code.

\subsection{PUFEM, refinement, and bit complexity}

The analytic inequalities in Sections~\ref{sec:reconstruction}--\ref{sec:descent} belong to the
partition-of-unity/GFEM tradition~\cite{MelenkBabuska96,BabuskaMelenk}.  The $h^{-1}$ derivative
factor and the $h^{r-1}$ global $H^1$ rate are not new.  The shared-DAG representation of hierarchical
refinement history used here is an explicit encoding choice, not a separate novelty claim.  The new
theorem is scoped to the combination: supplied finite local analytic witnesses are transported through
exact rational PUFEM synthesis and geometric refinement to a represented limit, while a self-contained
proof genealogy is retained under an explicit size/checking/lookahead bound.

Bit complexity of the approximant itself must be distinguished from proof overhead.  In a recent preprint, Mao and Xu
analyze bit-complexity rates across neural-network and classical approximation schemes, explicitly
including finite elements, and use metric-entropy comparisons to expose the role of encoding precision
~\cite{MaoXu2026}.  Theorem~\ref{thm:orderneutral} addresses a different comparison: relative to the
declared coefficient-wise $B_n$, the retained proof genealogy is $O(B_n)$ at the precision-relevant level.

\section{Version history}\label{sec:versions}

This manuscript is version~3 of arXiv:2506.22693 and supersedes versions~1--2 under the
same identifier.  The earlier programme already joined finite approximation with checkable
construction data, local-to-global assembly, categorical organization, and contextual
certified choice.  Version~3 retains that programme while separating claims that had been
run together.  Corrections are stated where the mathematical status changed, without
redescribing the earlier papers as claims about ambient classical non-existence.

Four revisions are unambiguously substantive.  The probes--models adjunction is withdrawn
as a theorem in its stated form: the earlier $F(T)$ depended on a target space not determined
by $T$, and the claimed natural hom-set bijection was not established.  The earlier
finite-rank projection formulation of UELAT is not retained as a universal approximation-
property statement.  The specific Chebyshev/B-spline numerical calculations identified as
incorrect during revision are omitted.  Priority language attached to generic certificate
stability is also not maintained; standard effective completion and strict-slack
certification are treated as background.  Universal gluing, by contrast, remains the
programme of explicit certified local-to-global construction, now realized by separate
reconstruction, defect, and refinement theorems under stated hypotheses.  CCP remains an internal discipline of admissible certified construction; version~3
expresses that discipline by generated closure and makes its internal scope explicit.

\renewcommand{\arraystretch}{1.15}
\begin{longtable}{>{\raggedright\arraybackslash}p{0.30\textwidth} >{\raggedright\arraybackslash}p{0.60\textwidth}}
\toprule
Versions 1--2 & Version 3 \\
\midrule
\endfirsthead
\toprule
Versions 1--2 & Version 3 \\
\midrule
\endhead
\midrule
\endfoot
\bottomrule
\endlastfoot
A stated probes--models adjunction organizing finite questions and analytic realizations & The adjunction is withdrawn at that level of generality; the intended architecture is retained through represented information, finite-code realization, chosen evidence lifts, the Grothendieck construction, and exact extensional sheaf semantics.\\
Proof-carrying analytic objects with analytic identity and construction genealogy treated as tightly integrated & Represented analytic identity and finite proof evidence are formally separated; selected genealogy remains proof-relevant structure above the represented object.\\
Uniform certificate stability and genealogy through limits, with a priority claim & Effective completion and generic strict-slack certification are treated as established; the new result is quantitative transport and retention of selected genealogy under explicit refinement and resource hypotheses.\\
Universal gluing as the certified local-to-global programme & The programme is retained; reconstruction, approximate gluing, and refinement are stated as distinct theorems with declared hypotheses, with the rational PUFEM compiler as the principal concrete realization.\\
CCP expressed in internal existence/admissibility language for a chosen topos or constructive universe & The internal scope is made explicit as generated genealogical closure under declared admissible constructors, with no inference about ambient classical non-existence.\\
Pathology exclusion in CCP-governed settings & Theorem~\ref{thm:preservation} gives invariant preservation for a specified generated universe; Corollary~\ref{cor:BT} is the concrete finite-measure exclusion example.\\
Universal embedding/linear approximation stated using finite-rank projection language & Effective linear isometric embedding into standard computable $C([0,1])$, followed by finite linear representation of embedded points; no universal approximation-property or bounded-approximation-property claim.\\
Specific Chebyshev/B-spline numerical worked material & The calculations identified as incorrect during revision are omitted; no current theorem depends on them.\\
\end{longtable}
\renewcommand{\arraystretch}{1.0}

\section{Formalization status and limitations}\label{sec:formal}

A dedicated public Rocq/Coq repository is maintained at \url{https://github.com/ipsissima/UELAT}.
Its \href{https://github.com/ipsissima/UELAT/blob/main/docs/FORMALIZATION_STATUS.md}{\texttt{docs/FORMALIZATION\_STATUS.md}}
is the authoritative correspondence table.  Only a paper result explicitly marked
\texttt{CHECKED-EXACT} there, for this manuscript's title/theorem numbering and at a
recorded audited commit, should be described as machine-checked.  Repository-wide CI,
legacy modules, or a row targeting a different manuscript snapshot are not evidence that a
current theorem has been formalized.  A checked claim additionally requires a reproducible
pinned toolchain, a formal statement at the same strength as the manuscript, inclusion in
the checked build, and an assumptions audit at the cited commit.  This manuscript makes no
blanket machine-checked claim.

Several limitations are intentional.
\begin{enumerate}[label=(\arabic*)]
\item Certificate sizes and verification times are encoding- and machine-relative.  No
Kolmogorov-style invariance theorem is asserted.
\item The proved compiler theorem is one-dimensional and uses a rational
piecewise-polynomial $W^{1,2}$ presentation, quasi-uniform dyadic refinement, bounded
overlap, and a fixed shared-DAG proof language.  Its geometric-summability mechanism is
dimension-independent: in an analogous $d$-dimensional quasi-uniform hierarchy with
$M_n\asymp h_n^{-d}$ and $h_n=h_0 2^{-n}$, one still has
$\sum_{j\le n}M_j=O(M_n)$.  Under analogous local $h^{r-1}$ approximation and a linear
bit schedule, the same bookkeeping would give the scale
$\varepsilon^{-d/(r-1)}(1+\log(1/\varepsilon))$.  Establishing a theorem in that setting
requires the higher-dimensional local approximation, shape-regularity, incidence, and
compiler hypotheses; adaptive variants require their own cumulative-complexity analysis,
as in the adaptive FEM literature~\cite{BinevDahmenDeVore,Stevenson}.
\item The theorem gives an upper bound relative to the conventional finite encoding; it
does not establish that the genealogy-only overhead is information-theoretically
necessary.
\item Qualitative local non-liftability is not claimed for the present Lipschitz evidence
grammar.  Theorem~\ref{thm:saturation} shows the opposite under its hypotheses.
\item The Grothendieck construction and exact sheaf placement are standard categorical
machinery.  No probes--models adjunction, enriched sheaf universal property, or stack
semantics is inferred from them.
\item The effective linear universality theorem preserves linear structure, but no
priority claim is needed for the paper's mathematical validity.  Metric computable
universality of $C([0,1])$ is established prior art.
\end{enumerate}

\section{Conclusion}\label{sec:conclusion}

Proof-carrying analytic approximation separates represented analytic identity from the
finite genealogy of a particular construction.  Generic complete strict-slack evidence
adds no new extensional carrier once full metric names are available, but selected
transformers, source use, proof size, and verification history remain intensional data.
The universal representation theorem shows that a common finite linear code language is
available across real computable Banach presentations, while the concrete Sobolev compiler
shows what is required to transport genealogy with quantitative control.

In the rational $W^{1,2}(0,1)$ presentation, supplied local approximation and overlap
witnesses are compiled through exact partition-of-unity synthesis rather than discarded
and regenerated from a semantic target.  The refinement construction therefore has
$Q_{\rm target}=0$.  More importantly, carrying the complete history up to the required
precision costs the same asymptotic bit order as the finest level alone:
geometric refinement makes all coarser $M_j\beta_j$ contributions sum to a constant
multiple of $M_n\beta_n$.  Proposition~\ref{prop:H6compiler} connects that statement to
the concrete one-cover compiler under transparent standard local-encoding assumptions;
Theorem~\ref{thm:orderneutral} keeps the formulation abstract enough to accommodate other
encodings.

Contextual Choice records the corresponding admissibility semantics.  Declared primitives,
admitted evidence-transforming constructors, and proved gluing/refinement interfaces
generate a certified genealogy whose preserved invariants follow by structural induction.
Exact sheaf semantics describes the extensional gluing boundary, while the proof-carrying
layer remembers data erased by that quotient.  The resulting framework is therefore
unified rather than foundationally maximal: representation, reconstruction, evidence
transport, refinement, and generated admissibility are distinct layers of one problem.
The version history in Section~\ref{sec:versions} records continuity and genuine revision;
none of the current theorems depends on that history.


\begin{thebibliography}{99}
\footnotesize
\setlength{\itemsep}{1.5pt}
\setlength{\parsep}{0pt}

\bibitem{AlbiacKalton}
F.~Albiac and N.~J.~Kalton,
\emph{Topics in Banach Space Theory}, 2nd ed., Graduate Texts in Mathematics 233,
Springer, 2016. \href{https://doi.org/10.1007/978-3-319-31557-7}{doi}.

\bibitem{BabuskaMelenk}
I.~Babu\v{s}ka and J.~M.~Melenk,
The partition of unity method,
\emph{International Journal for Numerical Methods in Engineering} 40 (1997), 727--758.
\href{https://doi.org/10.1002/(SICI)1097-0207(19970228)40:4<727::AID-NME86>3.0.CO;2-N}{doi}.

\bibitem{BagavievEtAl}
R.~Bagaviev, I.~I.~Batyrshin, N.~Bazhenov, D.~Bushtets, M.~Dorzhieva,
H.~T.~Koh, R.~Kornev, A.~G.~Melnikov, and K.~M.~Ng,
Computably and punctually universal spaces,
\emph{Annals of Pure and Applied Logic} 176 (2025), no.~1, 103491.
\href{https://doi.org/10.1016/j.apal.2024.103491}{doi}.

\bibitem{BanachTarski1924}
S.~Banach and A.~Tarski,
Sur la d\'ecomposition des ensembles de points en parties respectivement congruentes,
\emph{Fundamenta Mathematicae} 6 (1924), 244--277.

\bibitem{Dandelion}
H.~Becker, M.~Tekriwal, E.~Darulova, A.~Volkova, and J.-B.~Jeannin,
Dandelion: Certified approximations of elementary functions,
in \emph{13th International Conference on Interactive Theorem Proving (ITP 2022)},
LIPIcs 237, Article 6, 6:1--6:19, 2022.
\href{https://doi.org/10.4230/LIPIcs.ITP.2022.6}{doi}.

\bibitem{BinevDahmenDeVore}
P.~Binev, W.~Dahmen, and R.~DeVore,
Adaptive finite element methods with convergence rates,
\emph{Numerische Mathematik} 97 (2004), no.~2, 219--268.
\href{https://doi.org/10.1007/s00211-003-0492-7}{doi}.

\bibitem{BishopBridges}
E.~Bishop and D.~Bridges,
\emph{Constructive Analysis},
Grundlehren der mathematischen Wissenschaften 279, Springer, 1985.
\href{https://doi.org/10.1007/978-3-642-61667-9}{doi}.

\bibitem{BoldoEtAl2026}
S.~Boldo, F.~Cl\'ement, V.~Martin, M.~Mayero, and H.~Mouhcine,
A Rocq formalization of simplicial Lagrange finite elements,
\href{https://arxiv.org/abs/2604.20345v1}{arXiv:2604.20345v1 [cs.LO]}, 2026.

\bibitem{Brattka2005}
V.~Brattka,
On the Borel complexity of Hahn--Banach extensions,
\emph{Electronic Notes in Theoretical Computer Science} 120 (2005), 3--16.
\href{https://doi.org/10.1016/j.entcs.2004.07.011}{doi}.

\bibitem{Brattka2008}
V.~Brattka,
Borel complexity and computability of the Hahn--Banach theorem,
\emph{Archive for Mathematical Logic} 46 (2008), 547--564.
\href{https://doi.org/10.1007/s00153-007-0057-z}{doi}.

\bibitem{BrattkaDeBrechtPauly}
V.~Brattka, M.~de Brecht, and A.~Pauly,
Closed choice and a uniform low basis theorem,
\emph{Annals of Pure and Applied Logic} 163 (2012), 986--1008.
\href{https://doi.org/10.1016/j.apal.2011.12.020}{doi}.

\bibitem{BrattkaHertling}
V.~Brattka, P.~Hertling, and K.~Weihrauch,
A tutorial on computable analysis,
in \emph{New Computational Paradigms}, Springer, 2008, 425--491.
\href{https://doi.org/10.1007/978-0-387-68546-5_18}{doi}.

\bibitem{BrattkaSorg2026}
V.~Brattka and C.~Sorg,
Computability of the Hahn--Banach theorem revisited,
\href{https://arxiv.org/abs/2603.16802v1}{arXiv:2603.16802v1 [math.LO]}, 2026.

\bibitem{BrehardMahboubiPous}
F.~Br\'ehard, A.~Mahboubi, and D.~Pous,
A certificate-based approach to formally verified approximations,
in \emph{10th International Conference on Interactive Theorem Proving (ITP 2019)},
LIPIcs 141, Article 8, 8:1--8:19, 2019.
\href{https://doi.org/10.4230/LIPIcs.ITP.2019.8}{doi}.

\bibitem{ChihaniMillerRenaud}
Z.~Chihani, D.~Miller, and F.~Renaud,
A semantic framework for proof evidence,
\emph{Journal of Automated Reasoning} 59 (2017), no.~3, 287--330.
\href{https://doi.org/10.1007/s10817-016-9380-6}{doi}.

\bibitem{ColbrookHansen}
M.~J.~Colbrook and A.~C.~Hansen,
The foundations of spectral computations via the Solvability Complexity Index hierarchy,
\emph{Journal of the European Mathematical Society} 25 (2023), 4639--4718.
\href{https://doi.org/10.4171/JEMS/1289}{doi}.

\bibitem{DalLagoHofmann}
U.~Dal Lago and M.~Hofmann,
Realizability models and implicit complexity,
\emph{Theoretical Computer Science} 412 (2011), no.~20, 2029--2047.
\href{https://doi.org/10.1016/j.tcs.2010.12.025}{doi}.

\bibitem{DowneyMelnikov}
R.~G.~Downey and A.~G.~Melnikov,
Computably compact metric spaces,
\emph{Bulletin of Symbolic Logic} 29 (2023), no.~2, 170--263.
\href{https://doi.org/10.1017/bsl.2023.16}{doi}.

\bibitem{EdmundsTriebel}
D.~E.~Edmunds and H.~Triebel,
\emph{Function Spaces, Entropy Numbers, Differential Operators},
Cambridge Tracts in Mathematics 120, Cambridge University Press, 1996.
\href{https://doi.org/10.1017/CBO9780511662201}{doi}.

\bibitem{Enflo}
P.~Enflo,
A counterexample to the approximation problem in Banach spaces,
\emph{Acta Mathematica} 130 (1973), 309--317.
\href{https://doi.org/10.1007/BF02392270}{doi}.

\bibitem{ErnStephansenVohralik}
A.~Ern, A.~F.~Stephansen, and M.~Vohral\'ik,
Guaranteed and robust discontinuous Galerkin a posteriori error estimates for
convection--diffusion--reaction problems,
\emph{Journal of Computational and Applied Mathematics} 234 (2010), no.~1, 114--130.
\href{https://doi.org/10.1016/j.cam.2009.12.009}{doi}.

\bibitem{FourmanScott}
M.~P.~Fourman and D.~S.~Scott,
Sheaves and logic,
in \emph{Applications of Sheaves}, Lecture Notes in Mathematics 753, Springer, 1979,
302--401. \href{https://doi.org/10.1007/BFb0061824}{doi}.

\bibitem{HesthavenRozzaStamm}
J.~S.~Hesthaven, G.~Rozza, and B.~Stamm,
\emph{Certified Reduced Basis Methods for Parametrized Partial Differential Equations},
SpringerBriefs in Mathematics, Springer, 2016.
\href{https://doi.org/10.1007/978-3-319-22470-1}{doi}.

\bibitem{HofmannLedent}
M.~Hofmann and J.~Ledent,
A quantitative model for simply typed $\lambda$-calculus,
\emph{Mathematical Structures in Computer Science} 32 (2022), 777--793.
\href{https://doi.org/10.1017/S0960129521000256}{doi}.

\bibitem{Johnstone}
P.~T.~Johnstone,
\emph{Sketches of an Elephant: A Topos Theory Compendium},
2 vols., Oxford Logic Guides 43--44, Oxford University Press, 2002.

\bibitem{Kohlenbach}
U.~Kohlenbach,
\emph{Applied Proof Theory: Proof Interpretations and their Use in Mathematics},
Springer Monographs in Mathematics, Springer, 2008.
\href{https://doi.org/10.1007/978-3-540-77533-1}{doi}.

\bibitem{Lawvere73}
F.~W.~Lawvere,
Metric spaces, generalized logic, and closed categories,
\emph{Rendiconti del Seminario Matematico e Fisico di Milano} 43 (1973), 135--166.
\href{https://doi.org/10.1007/BF02924844}{doi}.

\bibitem{MaoXu2026}
T.~Mao and J.~Xu,
Do neural networks really beat the curse of dimensionality? A bit-complexity view,
\href{https://arxiv.org/abs/2608.01357v1}{arXiv:2608.01357v1 [cs.LG]}, 2026.

\bibitem{MelenkBabuska96}
J.~M.~Melenk and I.~Babu\v{s}ka,
The partition of unity finite element method: basic theory and applications,
\emph{Computer Methods in Applied Mechanics and Engineering} 139 (1996), 289--314.
\href{https://doi.org/10.1016/S0045-7825(96)01087-0}{doi}.

\bibitem{NeculaPCC}
G.~C.~Necula,
Proof-carrying code,
in \emph{Proceedings of the 24th ACM SIGPLAN--SIGACT Symposium on Principles of Programming Languages (POPL '97)},
ACM, 1997, 106--119.
\href{https://doi.org/10.1145/263699.263712}{doi}.

\bibitem{PourElRichards}
M.~B.~Pour-El and J.~I.~Richards,
\emph{Computability in Analysis and Physics}, Springer, 1989.

\bibitem{Sorg2026}
C.~Sorg,
Foundational analysis of the Solvability Complexity Index: the Weihrauch--SCI intermediate
hierarchy,
\href{https://arxiv.org/abs/2603.18955v3}{arXiv:2603.18955v3 [math.LO]}, 2026.

\bibitem{Stevenson}
R.~P.~Stevenson,
Optimality of a standard adaptive finite element method,
\emph{Foundations of Computational Mathematics} 7 (2007), no.~2, 245--269.
\href{https://doi.org/10.1007/s10208-005-0183-0}{doi}.

\bibitem{TWW}
J.~F.~Traub, G.~W.~Wasilkowski, and H.~Wo\'zniakowski,
\emph{Information-Based Complexity}, Academic Press, 1988.

\bibitem{vanOosten}
J.~van Oosten,
\emph{Realizability: An Introduction to its Categorical Side},
Studies in Logic and the Foundations of Mathematics 152, Elsevier, 2008.

\bibitem{Weihrauch}
K.~Weihrauch,
\emph{Computable Analysis: An Introduction},
Texts in Theoretical Computer Science, Springer, 2000.
\href{https://doi.org/10.1007/978-3-642-56999-9}{doi}.

\end{thebibliography}
\end{document}